\documentclass[mathpazo]{csiam-am}
\usepackage{graphicx}
\usepackage{subcaption}
\usepackage{amssymb}
\usepackage{mathtools}
\graphicspath{{Fig/}}
\usepackage{float}    
\usepackage{array}    
\usepackage{amsmath}
\usepackage{color}
\usepackage{tikz}
\usepackage{booktabs}
\usepackage{multirow}
\usepackage{caption}
\usepackage{graphicx}
\renewcommand\thisnumber{x}
\renewcommand\thisyear {20xx}
\renewcommand\thismonth{xx}
\renewcommand\thisvolume{x}

\def\norm#1{\Vert\,#1\,\Vert}
\def\Norm#1#2{\left\|\,#1\,\right\|_{#2}}

\def\no{{\nonumber}}

\newcommand{\tr}{\operatorname{tr}}

\def\ba{\mathbf{a}}

\newcommand{\ilQ}[1]{\norm{#1 }_{L}^2}
\newcommand{\ilQa}[1]{\norm{#1 }_{L_1}^2}

\newcommand{\iaQ}[1]{\norm{#1 }^2_{\psi_{L}}}
\newcommand{\ibQ}[1]{\norm{#1 }^2_{\psi_{L}^1}}

 \def\ba{\begin{equation}\begin{aligned}}
 \def\ed{\end{aligned}\end{equation}}

\def\mq{\psi}

\def\ba{\begin{align}}
\def\da{\end{align}}

\begin{document}
\title{Maximum bound principle and defect dynamics for Q-tensor gradient flows with low regularity integrators}

\author[O.~Author]{Only Author\corrauth}
\address{School of Mathematical Sciences, Beijing Normal University,
Beijing 100875, P.R. China}

\author[W. Hu and G. Ji]{Wenshuai Hu\affil{1}~and Guanghua Ji\affil{1}\comma\corrauth}
\address{\affilnum{1}\ Laboratory of Mathematics and Complex Systems, Ministry of Education and School of Mathematical Sciences,
         Beijing Normal University, Beijing 100875, China.}
\emails{{\tt 202431130052@mail.bnu.edu.cn} (W.~Hu),
         {\tt ghji@bnu.edu.cn} (G.~Ji)}

%
%

\begin{abstract}
The Landau-de Gennes (LdG) theory is a widely used thermodynamic continuum framework for modeling the behavior of ordered states and defects in liquid crystals with  a tensor-order parameter $Q$.
In this study, we develop and analyze first- and second-order low-regularity integrator (LRI) schemes  for  the $Q$-tensor gradient flow and prove the maximum bound principle.
In particular, through the reformulation of the LRI schemes, we establish rigorous modified energy dissipation laws for the LRI1a and LRI1b schemes, thereby filling a significant theoretical gap in the existing literature on LRI methods.
Moreover, this reformulation establishes a structural bridge between the LRI schemes and backward differentiation formula (BDF) methods, which opens up new possibilities for the construction and analysis of LRI-type  methods.
We then establish first- and second-order temporal convergence under $H^1$ and $H^2$ regularity assumptions, respectively.
 Several numerical experiments are presented  to validate our theoretical results and to simulate the evolution of defect dynamics.
\end{abstract}


\maketitle


\section{Introduction}

Liquid crystals (LCs) are a fascinating state of matter characterized by long-range orientational order. However, this order is rarely perfect throughout a sample, often manifesting as singularities known as topological defects. These defects are localized regions where the molecular alignment breaks down and can be categorized into point, line, or wall defects based on their topological characteristics.  Understanding and accurately modeling the dynamics of such defects are critical for both fundamental studies and practical applications \cite{deGennes1974, 10.1007/s10915-016-0177-x,wang2021modelling}.
The Landau-de Gennes  theory is a cornerstone thermodynamic continuum framework for describing ordered states and defects in liquid crystals. Specifically, it captures both the degree of molecular alignment and the preferred orientation of the molecules by employing a tensor order parameter $Q$ which is a symmetric, traceless second-rank tensor \cite{deGennes1974,mottramintroduction,Selinger2016}.

In this paper, we consider the Landau-de Gennes free energy functional for  nematic LCs in \cite{wang2021modelling,deGennes1974}, which is given by
\begin{alignat}{2}
	E(Q)=\int_\Omega \left[ \frac{c}{2}\|\nabla Q\|_F^2+\frac{\alpha}{2} \text{trace}(Q^2)-\frac{\beta}{3} \text{trace}(Q^3)+\frac{\gamma}{4} \text{trace}(Q^2)^2\right] d \Omega\label{free_energy},
\end{alignat}
where $\|\nabla Q\|_F^2 = \sum_{i,j,k=1}^d (\partial_k Q_{ij})^2$. Here, $\Omega \subset \mathbb{R}^d$ is a bounded domain, and we consider the physically relevant cases of $d=2$ and $d=3$.
 The parameter $\alpha = a \Delta \theta$, where $a > 0$, $\Delta \theta$ is the temperature difference relative to the transition temperature.
The parameters $c>0,\beta>0$ and $\gamma>0$ are material constants that characterize the elastic properties of the examined LCs \cite{nguyen2013refined,mottramintroduction}.

 From \cite{mottramintroduction,wang2021modelling}, the gradient flow equation of the free energy \eqref{free_energy} subject to the periodic boundary condition and the initial condition is given by
\begin{equation}\label{1.6a}
    \begin{aligned}
        	&Q_t = c\Delta Q-\alpha Q+\beta \left(Q^2-\frac{1}{d}\mathrm{tr}(Q^2)I\right)-\gamma \mathrm{tr}(Q^2)Q,  \\[5pt]
	&Q(0,\boldsymbol{x}) = Q_0(\boldsymbol{x}), \qquad~~~~ \boldsymbol{x} \in \Omega, \\[5pt]
	&Q(t,\cdot) \text{ is } \Omega\text{-periodic}, \qquad t \in [0,T].
    \end{aligned}
\end{equation}
where the term $-\frac{1}{d}\beta \operatorname{tr}(Q^2)I$ ensures the traceless property of the evolution.

 Extensive studies \cite{MariusPaicu2011SIAMJMA,  Francosco2014SIAMJ.MA, MariusPaicu2012Arch,  HelmutAbels2016AdvancesInDE} have been conducted on analyzing the regularity  of solutions and the dissipation of energy to the coupled Navier-Stokes and Q-tensor model. Specifically, the work in \cite{weak_solutions_Q} not only establishes the existence of weak solutions for the coupled model but also proves that the Q-tensor satisfies the maximum bound preserving (MBP) property, provided that the condition $b \leq a^2$ holds, where $b = \frac{\beta^2}{\gamma^2}-\frac{2\alpha}{\gamma}$ and $a^2$ is the supremum of the Frobenius norm of the initial data $Q_0$.
Moreover, for nematic LCs, the work in \cite{hu2016disclination,majumdar2010equilibrium} shows that the eigenvalues of the Q-tensor are physically constrained to the intervals $[-1/3, 2/3]$ and $[-1/2, 1/2]$ for three-dimensional (3D) and two-dimensional (2D) systems, respectively.

Numerical schemes for calculating the Q-tensor gradient flows   and  the coupled system   have also attracted significant attention in recent years.
Hou et al \cite{qiaozhonghua}  proposed first- and second-order  stabilized exponential scalar auxiliary variable $\rm (SAV)$ schemes for the Q-tensor flows of
LCs, which  preserve the MBP  and energy dissipation.
Zhao et al \cite{10.1007/s10915-016-0177-x} proposed a linear second-order energy stability scheme for the Navier-Stokes coupled Q-tensor model based on the IEQ method.
Ji et al \cite{huang2024errorestimateorderenergy}\cite{10.4208/eajam.291018.070619} proposed a SAV type energy stability scheme and provided an error estimate for the first order scheme.

The boundedness of the Q-tensor's eigenvalues ensures that its Frobenius norm also satisfies the MBP property.  Furthermore, the presence of defects in typical initial configurations leads to limited  solution regularity.  Consequently,  numerical schemes should preserve the MBP and maintain their accuracy under low-regularity condition.
  Many numerical schemes are available for solving scalar parabolic equations that can preserve the MBP, such as the ETD schemes \cite{Ju2018,Du2019,du2021,Li2020},  the IFRK schemes \cite{Integrating2022,Li2021,Li2023}, and the implicit-explicit  schemes \cite{wang2021,TangYang2016,hou2023linear}.   Specially, the authors of \cite{Li2021} mentioned that the ETD schemes could be extended to the the matrix-valued Allen-Cahn equation and Liu et al  \cite{liu2024maximum} provided the detailed proofs of the unconditional MBP and the energy dissipation. However,  these schemes rely on high solution regularity.

Recently, low regularity integrators (LRIs) have emerged as a powerful framework for addressing problems where the solutions lack sufficient smoothness \cite{MR4613580,SchratzWangZhao2021,Rousset2021}. By rigorously preserving the MBP and energy stability under minimal regularity assumptions, LRIs are particularly well-suited for simulating systems with inherent singularities, such as defects in liquid crystals \cite{Doan2022LowRI,OstermannSu2019}.
While LRI schemes have been successfully applied to gradient flows in recent years, the rigorous energy dissipation analysis is still lacking, which is a fundamental property for gradient flow models \cite{FengMaierhoferSchratz2024,Doan2022LowRI}.
In this paper, we reformulate the LRI schemes as a  backward differentiation formula (BDF) and establish rigorous modified energy dissipation laws for the first-order LRI schemes, thereby filling this significant theoretical gap.
 Moreover, this reformulation bridges the gap between LRI schemes and BDF methods, which opens up new possibilities for the construction and analysis of LRI-type methods.


The main contributions are summarized as follows:
\begin{itemize}
	\item We construct first- and second-order low-regularity integrators for the Q-tensor gradient flow system \eqref{1.6a} and rigorously prove that these schemes preserve the MBP.
	\item We reformulate the  LRI schemes  and establish rigorous modified energy dissipation laws for the first-order LRI schemes, thereby filling a significant theoretical gap in the existing literature on LRI methods.
	\item We rigorously establish first-order convergence under an $H^1$ assumption and second-order convergence under an $H^2$ assumption, and numerically demonstrate the enhanced robustness of our method compared to other classical exponential integrator schemes.
\end{itemize}

The rest of this paper is organized as follows. In Section \ref{section2}, we present four $\rm LRI$ schemes for the  Q-tensor model. In Section \ref{section3}, we rigorously analyze these schemes in terms of the  ${\rm MBP}$ and energy stability. In Section \ref{section4},  we provide detailed proofs and derive error estimates for the proposed schemes.
In Section \ref{section6}, we present numerical experiments to validate the theoretical findings and demonstrate the robustness of the schemes.  We then conduct further numerical simulations to simulate the defect dynamics of  Q-tensor gradient flows.  Finally, in section \ref{section5}, we draw  conclusions.

\section{Numerical schemes}\label{section2}~\\
For the second-order tensor $A,B$ and their gradients   $\nabla A$ and $\nabla B$,
 the Frobenius inner product $(\cdot,\cdot)_F$ and the Frobenius norm $\|\cdot\|_F$  are respectively defined as
\begin{align*}
    (A,B)_F&=A : B=\sum_{i,j} A_{ij} B_{ij},\quad (  \nabla A,\nabla B)_F=\sum_{i,j,k} \partial_k A_{ij} \partial_k B_{ij},\\
	\|A\|_F^{2}&=(A,A)_F, \quad\|\nabla A\|_F^{2}=(\nabla A,\nabla A)_F.
\end{align*}

Based on the symmetric and traceless properties of the Q-tensor, the   Q-tensor space $\mathbb{R}_{s,0}^{d \times d}$ is defined as
\begin{equation}
	\mathbb{R}_{s,0}^{d \times d} = \left\{Q \in \mathbb{R}^{d \times d}\left\lvert\,Q_{i j}=Q_{j i},  Q_{i i}=0 \right.\right\}.
\end{equation}
For a bounded domain $\Omega \subset \mathbb{R}^{d}$, we define $\mathcal{Z}=L^{\infty}\left(\bar{\Omega} ; \mathbb{R}_{s,0}^{d \times d}\right)$ as the space of essentially bounded functions from $\bar{\Omega}$ to   $\mathbb{R}_{s,0}^{d \times d}$ equipped with the norm
\begin{align*}
	\|Q\|_\mathcal{Z}&= \max _{\boldsymbol{x} \in \bar{\Omega}}\|Q(\boldsymbol{x})\|_F, \quad \forall Q \in \mathcal{Z}.
\end{align*}
For $Q\in \mathcal{Z}$,  the uniform spectral norm  $\|\cdot\|_{\infty,2}$  is defined as
\begin{equation}
    \|Q\|_{\infty,2} = \max_{\boldsymbol{x} \in \bar{\Omega}} \|Q(\boldsymbol{x})\|_2,
\end{equation}
where $\|\cdot\|_2$ represents the standard spectral norm (matrix 2-norm) at a specific point $\boldsymbol{x}$.
The space $\mathcal{X}=L^{2}\left(\bar{\Omega} ; \mathbb{R}_{s,0}^{d \times d}\right)$ is defined as the Sobolev space equipped with the norm
\begin{align*}
\|Q\|_\mathcal{X}=(Q,Q)_\mathcal{X}^{\frac{1}{2}}\coloneqq \int_{\Omega}\|Q(\boldsymbol{x})\|_F^{2} d\boldsymbol{x}^{\frac{1}{2}}, \quad \forall Q \in \mathcal{X},
\end{align*}
and $H^1\coloneqq  W^{1,2}\left(\bar{\Omega} ; \mathbb{R}_{s,0}^{d \times d}\right),H^2\coloneqq  W^{2,2}\left(\bar{\Omega} ; \mathbb{R}_{s,0}^{d \times d}\right)$.

In the following sections, for ease of notation, we denote the gradient flow equation \eqref{1.6a} as $Q_t = c\Delta Q + f(Q)$, where $f(Q)=-\alpha Q+\beta (Q^2-\frac{1}{d}\text{trace}(Q^2)I)-\gamma \text{trace}(Q^2)Q$.  When $d=2$,  the second term in $f(Q)$ vanishes. Therefore, throughout the remainder of this manuscript, we let $d=3$, and all the results can be easily extended to  the case in which $d=2$ simply by  letting $\beta=0$.

Utilizing Duhamel's formula \cite{Rousset2021},   we can express the exact solution of \eqref{1.6a} as follows:
\begin{align}
	Q(t)=e^{ct \Delta} Q_0+\int_{0}^{t}e^{c(t-\xi) \Delta}f(Q(\xi))d\xi.\label{1.9}
\end{align}
Suppose the  terminal time is set to $T>0$, $\tau=\frac{T}{n}$, where $n$ is a positive integer.   We then define the time intervals as $t_m=m\tau$, $m=0,1,\ldots,n$ and the numerical solution at each time step as $Q_m$.
For an interval $[t_m, t_{m+1}]$, \eqref{1.9} can be written as
\begin{align}
	Q(t_{m+1})=e^{c\tau \Delta} Q(t_m)+\int_{0}^{\tau}e^{c(\tau-\xi) \Delta}f(Q(t_m+\xi))d\xi.\label{1.9a}
\end{align}

By using $e^{c(\tau-\xi) \Delta}f(e^{c\xi \Delta} Q)$ to approximate $e^{c(\tau-\xi) \Delta}f(Q(t_m+\xi))$, we obtain the following two first-order numerical schemes.
Taking $\xi=0$, the LRI1a scheme can be expressed as
\begin{align}
	Q_{m+1} = e^{c\tau \Delta} Q_m + \tau e^{c\tau \Delta} f(Q_m). \label{lri1a}
\end{align}
Taking $\xi=\tau$, the LRI1b scheme can be expressed as
\begin{align}
	Q_{m+1} = e^{c\tau \Delta} Q_m + \tau f(e^{c\tau \Delta} Q_m).\label{lri1b}
\end{align}
This scheme evaluates the nonlinear term $f(Q)$ at the intermediate state $e^{c\tau \Delta} Q_m$.

For the second-order LRI2a scheme,
using $e^{c\xi \Delta} Q_m+\xi f(Q_m)$ to approximate $Q(t_m+\xi)$
as in  \cite{Doan2022LowRI},  we can obtain
\begin{align}
	Q_{m+1} = e^{c\tau \Delta} Q_m + \frac{\tau}{2} \left[ e^{c\tau \Delta} f(Q_m) + f(e^{c\tau \Delta} Q_m) \right]
	+ \frac{\tau^2}{2} e^{c\tau \Delta} \frac{\partial f}{\partial Q}(Q_m) : f(Q_m), \label{lri2}
\end{align}
where $\frac{\partial f}{\partial Q}(Q_m)$ is the derivative of $f$ with respect to $Q$ evaluated at $Q_m$.

The other second-order LRI2b scheme is derived from the first-order LRI1b scheme.
Utilizing the  trapezoidal rule to approximate the second term on  the right-hand side of \eqref{1.9a}, we have
\begin{align*}
	Q_{m+1} = e^{c\tau \Delta} Q_m + \frac{\tau}{2} \left[ e^{c\tau \Delta} f(Q_m) + f(Q(t_m+\tau)) \right].
\end{align*}
By using LRI1b to approximate $Q(t_m+\tau)$ and Taylor expansion,  we can obtain the second-order LRI2b scheme:
\begin{align}
	Q_{m+1} = e^{c\tau \Delta} Q_m + \frac{\tau}{2} \left[ e^{c\tau \Delta} f(Q_m) + f(e^{c\tau \Delta} Q_m) \right]
	+ \frac{\tau^2}{2} \frac{\partial f}{\partial Q}(e^{c\tau \Delta}Q_m) : f(e^{c\tau \Delta}Q_m). \label{lri2b}
\end{align}

Through a series of tedious calculations, we have
\begin{align}
	\frac{\partial f}{\partial Q}(Q):f(Q)&=(-\alpha- \gamma \text{trace}(Q^2))f(Q)-2 \gamma(f(Q):Q)Q +2\beta( f(Q)Q-\frac{1}{3}(f(Q):Q)I).
	\label{1.16}
\end{align}
Similarly, when  $\hat{Q}$ is used to denote $e^{c\tau \Delta} Q$,  for $\frac{\partial f}{\partial Q}(e^{c\tau \Delta}Q_m) : f(e^{c\tau \Delta}Q_m)$, we have that
\begin{align}
	\frac{\partial f}{\partial {Q}}(\hat{Q}):f(\hat{Q})=(-\alpha- \gamma \text{trace}(\hat{Q}^2))f(\hat{Q})-2 \gamma(f(\hat{Q}):\hat{Q})\hat{Q} +2\beta( f(\hat{Q})\hat{Q}-\frac{1}{3}(f(\hat{Q}):\hat{Q})I).
	\label{1.16a}
\end{align}
\begin{remark}
	In \eqref{1.9a}, when we  use $Q(t_m)$ to approximate $Q(t_m+\xi)$ and using the Runge-Kutta method, we can obtain the  ETD1 and ETDRK2 schemes as listed in \cite{Du2019,du2021}, but they needs the assumption that $Q$ are
	respectively $H^2$  and $H^4$ regularity in space.
	To overcome this limitation, we consider constructing LRI schemes that  only require $H^1$  and $H^2$ regularity of $Q$, making them applicable to the problems   involving singularities and defects.
	Moreover, our theoretical analysis can be easily extended to the ETD schemes.
\end{remark}
\section{Properties of the semidiscrete numerical solutions}\label{section3}
\subsection{Discrete maximum bound principle}~\\
In this section, we  prove that the proposed schemes satisfy the MBP  under certain conditions.
\begin{lemma} \label{1.17}
    For any $\lambda > 0$, the linear operator $\Delta$ generates a contraction semigroup $\{e^{t\Delta}\}_{t \geq 0}$ on the space $\mathcal{X}$ (resp. $\mathcal{Z}$) if the following resolvent estimates hold:
    \begin{align}
                \lambda \|W\|_{\mathcal{X}} &\leq \|(\lambda I-\Delta)W\|_{\mathcal{X}}, \quad \forall W \in \mathcal{X}. \label{1.17b}\\
        \lambda \|W\|_{\mathcal{Z}} &\leq \|(\lambda I-\Delta)W\|_{\mathcal{Z}}, \quad \forall W \in \mathcal{Z}, \label{1.17a}
    \end{align}
\end{lemma}
\begin{proof}
	For $W \in \mathcal{X}$ and  $W(\boldsymbol{x})=\{w_{ij}(\boldsymbol{x}),i=1,2,3,j=1,2,3 \}$, we have that
\begin{align*}
(W(\boldsymbol{x}),\Delta W(\boldsymbol{x}))_{\mathcal{X}}&=\sum_{i=1}^3 \sum_{j=1}^3\int_{\Omega}^{} w_{ij}(\boldsymbol{x})\Delta w_{ij}(\boldsymbol{x})d \Omega=\sum_{i=1}^3 \sum_{j=1}^3\int_{\Omega}^{} -|\nabla w_{ij}(\boldsymbol{x})|^2d \Omega\leq 0.
\end{align*}
Then for any $\lambda  > 0$, using the Cauchy$–$Schwarz inequality,  we have
\begin{align}
	\lambda\| W(\boldsymbol{x})\|_{\mathcal{X}}^2\leq &\lambda \sum_{i=1}^3 \sum_{j=1}^3\int_{\Omega}^{}w_{ij}^2(\boldsymbol{x})d \Omega-\sum_{i=1}^3 \sum_{j=1}^3\int_{\Omega}^{} w_{ij}(\boldsymbol{x})\Delta w_{ij}(\boldsymbol{x})d \Omega\\=&
	\sum_{i=1}^3 \sum_{j=1}^3\int_{\Omega}^{} w_{ij}(\boldsymbol{x})(\lambda I-\Delta)w_{ij}(\boldsymbol{x})d \Omega
	\\\leq&\left(\sum_{i=1}^3 \sum_{j=1}^3\int_{\Omega}^{} |w_{ij}(\boldsymbol{x})|^2 d \Omega\right)^{\frac{1}{2}} \left(\sum_{i=1}^3 \sum_{j=1}^3\int_{\Omega}^{} |(\lambda I-\Delta)w_{ij}(\boldsymbol{x})|^2 d \Omega\right)^{\frac{1}{2}}
	\\=& \| W(\boldsymbol{x})\|_{\mathcal{X}} \|(\lambda I-\Delta) W(\boldsymbol{x})\|_{\mathcal{X}}. \label{2.3}
\end{align}
From \eqref{2.3}, we have
\begin{align}
	\lambda \| W(\boldsymbol{x})\|_{\mathcal{X}}\leq  \|(\lambda I-\Delta) W(\boldsymbol{x})\|_{\mathcal{X}}. \label{3.3}
\end{align}
Based on the proof of \cite[Lemma 2.1]{du2021}, \eqref{3.3} satisfies  the property (ii). Then, we can deduce the
contraction property of the operator $e^{t  \Delta}$ on the space $\mathcal{X}$.

	For  $W \in \mathcal{Z}$, the inequality \eqref{1.17a} has been proven and readers can  refer to \cite[Lemma 2.3]{liu2024maximum} for the associated details.
\end{proof}
\begin{lemma}\label{1.1}
    Suppose we are given a fixed terminal time $T>0$, a positive integer $n$, and a time step size $\tau=\frac{T}{n}$. Choose an explicit upper bound $\eta > \frac{\beta + \sqrt{\beta^2 - 24\alpha\gamma}}{2\sqrt{6}\gamma}$.
    For any $Q \in \mathcal{Z}$ satisfying $\|Q\|_F \leq \eta$, when $\tau \leq \tau^*= \frac{1}{|\alpha| + \beta \eta + 3\gamma \eta^2}$, we have:
    \begin{align}
        \|Q+\tau f(Q)\|_F^2 \leq \eta^2. \label{3,9}
    \end{align}
    Moreover, when $\tau \leq \frac{\alpha - \frac{\beta}{\sqrt{6}}\eta + \gamma \eta^2}{|\alpha| + \frac{\beta}{\sqrt{6}}\eta + \gamma \eta^2}\tau^*$, we have:
    \begin{align}
        \|Q+\tau f(Q)+\tau^2\frac{\partial f}{\partial Q}(Q): f(Q)\|_F^{2} \leq \eta^2. \label{3,10}
    \end{align}
\end{lemma}
\begin{proof}
    By the definition of the Frobenius norm and the properties of traceless symmetric tensors $Q \in \mathcal{Z}$, we apply the Cayley-Hamilton theorem to derive the exact identity:
    \begin{align}
        \left\| Q^2 - \frac{1}{3}\tr(Q^2)I \right\|_F^2 = \frac{1}{6}\|Q\|_F^4. \label{eq_CH_identity}
    \end{align}
    Consider the term $\|Q+\tau f(Q)\|_F$. By substituting the explicit form of $f(Q)$ and applying the triangle inequality, we obtain:
    \begin{align}
        \|Q+\tau f(Q)\|_F &= \left\| (1 - \tau\alpha - \tau\gamma\|Q\|_F^2)Q + \tau\beta\left(Q^2 - \frac{1}{3}\tr(Q^2)I\right) \right\|_F \no \\
        &\leq |1 - \tau\alpha - \tau\gamma\|Q\|_F^2| \|Q\|_F + \tau\beta \left\| Q^2 - \frac{1}{3}\tr(Q^2)I \right\|_F. \label{eq_triangle}
    \end{align}
    For $\tau \leq \tau_0 $ and $\|Q\|_F \leq \eta$, we guarantee that the coefficient $1 - \tau\alpha - \tau\gamma\|Q\|_F^2 \geq 0$. Letting $x = \|Q\|_F \in [0, \eta]$, we introduce the scalar function $g(x)$:
    \begin{align}
        \|Q+\tau f(Q)\|_F \leq g(x) &= (1 - \tau\alpha - \tau\gamma x^2)x + \frac{\tau\beta}{\sqrt{6}}x^2 \no \\
        &= x + \tau\left(-\alpha + \frac{\beta}{\sqrt{6}}x - \gamma x^2\right)x \no \\
        &= x + \tau h(x) x.
    \end{align}

    To ensure that $g(x)$ is monotonically increasing on $[0, \eta]$, we examine its derivative $g'(x) = 1 + \tau\left(-\alpha + \frac{2\beta}{\sqrt{6}}x - 3\gamma x^2\right)$. Since $\tau \leq \tau^*$, it follows that $g'(x) \geq 0$, making $g(x)$ strictly increasing. By our threshold choice, we guarantee $h(\eta) < 0$ (which holds unconditionally for any $\eta$ if $\beta^2 \leq 24\alpha\gamma$, and explicitly requires $\eta > \frac{\beta + \sqrt{\beta^2 - 24\alpha\gamma}}{2\sqrt{6}\gamma}$ if $\beta^2 > 24\alpha\gamma$). Consequently, it directly follows that $g(\eta) < \eta$.
    Therefore, for all $x \in [0, \eta]$, we have
    \begin{align}
        \|Q+\tau f(Q)\|_F \leq g(x) \leq g(\eta) = \eta + \tau h(\eta) \eta \leq \eta, \quad \forall x \in [0, \eta]. \label{eq_g_bound},
    \end{align}
    which yields the conclusion \eqref{3,9}.

    Next, we consider the higher-order term. We first establish a bound for $\|f(Q)\|_F$ when $\|Q\|_F \leq \eta$:
    \begin{align}
        \|f(Q)\|_F \leq (\alpha + \gamma\|Q\|_F^2)\|Q\|_F + \frac{\beta}{\sqrt{6}}\|Q\|_F^2 \leq \eta\left(\alpha + \frac{\beta}{\sqrt{6}}\eta + \gamma \eta^2\right) = C_f \eta.
    \end{align}
    Using the explicit representation of $\frac{\partial f}{\partial Q}(Q):f(Q)$, we can derive the following bound:
    \begin{align}
        \left\|\frac{\partial f}{\partial Q}(Q):f(Q)\right\|_F &\leq \|(|\alpha| + \gamma\|Q\|_F^2)f(Q)\|_F + 2\gamma\|f(Q)\|_F\|Q\|_F^2 + 2\beta\|f(Q)Q\|_F \no \\
        &\leq (|\alpha| + \gamma \eta^2)\|f(Q)\|_F + 2\gamma \eta^2\|f(Q)\|_F + 2\beta \eta\|f(Q)\|_F \no \\
        &= (|\alpha| + 2\beta \eta + 3\gamma \eta^2)\|f(Q)\|_F \no \\
        &\leq (|\alpha| + 2\beta \eta + 3\gamma \eta^2) C_f \eta = C_{\partial} \eta.\label{lem2.2a}
    \end{align}
    Consequently, we can bound the entire higher-order term sequence:
    \begin{align}
        &\left\|Q+\tau f(Q)+\tau^2\frac{\partial f}{\partial Q}(Q): f(Q)\right\|_F \no \\
        \leq & \|Q+\tau f(Q)\|_F + \tau^2 \left\|\frac{\partial f}{\partial Q}(Q): f(Q)\right\|_F \no \\
        \leq & x + \tau h(x) x + \tau^2 C_{\partial} x \no \\
        = & x + \tau(h(x) + \tau C_{\partial})x =: g_1(x).
    \end{align}
    Since $g(x)$ and $\tau^2 C_{\partial} x$ are both monotonically increasing on $[0, \eta]$, their sum $g_1(x)$ is strictly monotonically increasing. Thus:
    \begin{align}
        \left\|Q+\tau f(Q)+\tau^2\frac{\partial f}{\partial Q}(Q): f(Q)\right\|_F \leq g_1(x) \leq g_1(\eta) = \eta + \tau(h(\eta) + \tau C_{\partial})\eta.
    \end{align}
    To ensure that $g_1(\eta)\leq \eta$, we require $h(\eta) + \tau C_{\partial} \leq 0 \implies \tau \leq \frac{-h(\eta)}{C_{\partial}}=\frac{-h(\eta)}{C_{f}}\tau_0$. Then we have
    \begin{align}
        \left\|Q+\tau f(Q)+\tau^2\frac{\partial f}{\partial Q}(Q): f(Q)\right\|_F \leq \eta.
    \end{align}
\end{proof}
\begin{remark}
	From  \eqref{3,9}, letting $\kappa =\frac{1}{\tau},\kappa\geq \frac{1}{\tau_0}$, we can obtain that $\|\kappa Q+ f(Q)\|_F^2\leq \kappa^2 \eta^2$ holds.  Based on the proof of \cite[Theorem 3.1]{liu2024maximum}, we can deduce that the ETD schemes  preserve the MBP unconditionally.
\end{remark}
\begin{theorem}\label{mbp}
For any $Q_0 \in \mathcal{Z}$ satisfying $\|Q_0\|^2_{\mathcal{Z}}\leq \eta^2$, when $\tau\leq\tau_0=\min\{\tau^*, \frac{-h(\eta)}{C_{f}}\tau^*\}$, the solution $Q_m$ generated by  schemes \eqref{lri1a}-\eqref{lri2b} satisfies
	\begin{align}
		\| Q_m\|_{\mathcal{Z}}^2\leq \eta^2, \qquad m=1,2,\ldots,n.
	\end{align}
\end{theorem}
\begin{proof}
	Suppose that $\|Q_{m-1}\|_F^2\leq \eta^2$, we  only need to show that $\|Q_m\|_F^2\leq \eta^2$ holds.
	Taking the supremum norm $\|\cdot\|_{\mathcal{Z}}$ on both sides of the LRI1a scheme \eqref{lri1a} and using Lemma \ref{1.17}, we obtain
	\begin{align}
		\|Q_{m}\|_\mathcal{Z}^{2}&=\|e^{c\tau \Delta} Q_{m-1}+\tau e^{c\tau \Delta}f(Q_{m-1})\|_\mathcal{Z}^{2}\no\\
		                 &\leq \|Q_{m-1}+\tau f(Q_{m-1})\|_\mathcal{Z}^{2}\leq \eta^2. \label{1.25c}
	\end{align}

	Taking  $\|\cdot\|_{\mathcal{Z}}$ on both sides of the LRI1b scheme \eqref{lri1b} and using Lemma \ref{1.17}, we obtain
	\begin{align}
		\|Q_{m}\|_\mathcal{Z}^{2}=\|e^{c\tau \Delta} Q_{m-1}+\tau f(e^{c\tau \Delta}Q_{m-1})\|_\mathcal{Z}^{2}\leq \eta^2.
	\end{align}

	Taking  $\|\cdot\|_{\mathcal{Z}}$ on both sides of the LRI2a scheme \eqref{lri2} and using Lemma \ref{1.17}, we obtain
	 \begin{align}
	    \|Q_{m}\|_\mathcal{Z}^{2}&=\|e^{c\tau \Delta} Q_{m-1}+\frac{1}{2}\tau[ e^{c\tau \Delta}f(Q_{m-1})+f(e^{c\tau \Delta}Q_{m-1})]+\frac{1}{2}\tau^2 e^{c\tau \Delta}\frac{\partial f}{\partial Q}(Q_{m-1}):f(Q_{m-1})\|_\mathcal{Z}^{2}\no\\
	    &\leq \frac{1}{2}\|e^{c\tau \Delta} Q_{m-1}+\tau f(e^{c\tau \Delta}Q_{m-1})\|_\mathcal{Z}^{2}\no \\
	    &\quad+\frac{1}{2}\|e^{c\tau \Delta} Q_{m-1}+\tau e^{c\tau \Delta}f(Q_{m-1})+\tau^2 e^{c\tau \Delta}\frac{\partial f}{\partial Q}(Q_{m-1}):f(Q_{m-1})\|_\mathcal{Z}^{2}\no \\
	    &\leq \frac{1}{2}\|Q_{m-1}\|_\mathcal{Z}^{2}+\frac{1}{2}\| Q_{m-1}+\tau f(Q_{m-1})+\tau^2\frac{\partial f}{\partial Q}(Q_{m-1}):f(Q_{m-1})\|_\mathcal{Z}^{2}\leq \eta^2.
	 \end{align}

	 Taking  $\|\cdot\|_{\mathcal{Z}}$ on both sides of the LRI2b scheme \eqref{lri2} and using Lemma \ref{1.17}, we obtain
	 \begin{align}
	    \|Q_{m}\|_\mathcal{Z}^{2}
	    &\leq \frac{1}{2}\|e^{c\tau \Delta} Q_{m-1}+\tau f(e^{c\tau \Delta}Q_{m-1})+\tau^2 \frac{\partial f}{\partial Q}(e^{c\tau \Delta}Q_{m-1}):f(e^{c\tau \Delta}Q_{m-1})\|_\mathcal{Z}^{2}\no \\
	    &\quad+\frac{1}{2}\|e^{c\tau \Delta} Q_{m-1}+\tau e^{c\tau \Delta}f(Q_{m-1})\|_\mathcal{Z}^{2}\no \\
	    &\leq \frac{1}{2}\|e^{c\tau \Delta}Q_{m}\|_\mathcal{Z}^{2}+\frac{1}{2}\| Q_{m-1}+\tau f(Q_{m-1})\|_\mathcal{Z}^{2}\leq \eta^2.\label{1.25a}
	 \end{align}
	  \eqref{1.25c}-\eqref{1.25a} verify that $\|Q_m\|_\mathcal{Z}^{2}\leq \eta^2$ holds for $m=1,2,\ldots,n$. The proof is complete.
\end{proof}
\begin{remark}\label{re3.1}
	From the definitions of the $\mathcal{Z}$ norm and the uniform spectral norm,  we can derive$\|Q\|_{\infty, 2}\leq \Norm{Q}{\mathcal{Z}}$,
which implies that any bound established for $Q$ in the $\mathcal{Z}$ norm (specifically the MBP) naturally extends to the uniform spectral norm.
 \end{remark}
 Since $f(Q),\frac{\partial f}{\partial Q}(Q):f(Q) $  are  polynomials in Q,    $f(Q)$ and $\frac{\partial f}{\partial Q}(Q):f(Q) $ are continuous functions of Q.  For any  $Q,Q_1 \in \mathcal{Z}$ satisfying $\|Q\|^2_F\leq \eta^2$, for any $\boldsymbol{x} \in \bar{\Omega}$, there exist constants $C_1,C_2,C_3$ that,   only depend on the coefficients of $f$ and $\eta$, where
	\begin{align}
		\|\nabla f(Q)\|_F&\leq C_2(\eta)\Norm{\nabla Q}{F},
\|\nabla \frac{\partial f}{\partial Q}(Q):f(Q)\|_F\leq C_3(\eta)\Norm{\nabla Q}{F},\label{lem2.3b}\\
&\|f(Q_1)-f(Q)\|_F\leq C(\eta)\|Q_1-Q\|_F,\label{lem2.2b}\\
&\|\frac{\partial f}{\partial Q}(Q_1)-\frac{\partial f}{\partial Q}(Q)\|_F\leq C_1(\eta)\|Q_1-Q\|_F.\label{lem2.3a}
	\end{align}
 \subsection{Semidiscrete energy stability}~\\
In this section, we prove the semidiscrete energy stability of  four  schemes \eqref{lri1a}-\eqref{lri2b}.

Let $L=-c \Delta$.
Acting $\frac{1}{\tau}(\mq(\tau L)+\tau L)$ on the both sides of \eqref{lri1a}, where $\mq(z) = z / (e^z - 1)$, we have
\begin{align}
    \frac{1}{\tau}(\mq(\tau L)+\tau L)(Q_{m+1}-e^{-\tau  L}Q_{m})&=(\mq(\tau  L)+\tau  L)e^{-\tau  L}f(Q_{m}). \label{lri1a_1}
\end{align}
For the left-hand side (LHS) of \eqref{lri1a_1}, we have
\begin{align}
    \quad \frac{1}{\tau}(\mq(\tau L)+\tau L)(Q_{m+1}-e^{-\tau L}Q_{m})= \frac{1}{\tau}\mq(L)(Q_{m+1}-Q_{m})+L Q_{m+1},
\end{align}
where we have used the identity $\mq(z)e^{z} = z + \mq(z)$.
For the right-hand side (RHS) of \eqref{lri1a_1}, we have
\begin{align}
    \frac{1}{\tau}(\mq(\tau L)+\tau L)e^{-\tau L}f(Q_{m})=\mq(\tau L)f(Q_{m}).
\end{align}
Then, we can derive the following equivalent form of the LRI1a scheme \eqref{lri1a}:
\begin{align}
    \frac{1}{\tau}\mq(\tau L)(Q_{m+1}-Q_{m})+L Q_{m+1}=\mq(\tau L)f(Q_{m}). \label{lri1a_2}
\end{align}
Dividing by $\mq(\tau L)$, we get
\begin{align}
    \frac{1}{\tau}(Q_{m+1}-Q_{m})+L_1 Q_{m+1}=f(Q_{m}), \label{lri1a_3}
\end{align}
where $L_1 =  \frac{L}{\mq(\tau L)}= \frac{e^{\tau L}-I}{\tau}$.
Similarly, we can derive the equivalent forms of the other scheme \eqref{lri1b}-\eqref{lri2b} as follows:
\begin{align}
    \frac{1}{\tau} \mq(\tau L)(Q_{m+1}-Q_{m}) +L Q_{m+1}&=\mq(\tau L)f(e^{-\tau L}Q_{m})+\tau Lf(e^{-\tau L}Q_{m}), \label{lri1b_2}\\
    \frac{1}{\tau} \mq(\tau L)(Q_{m+1}-Q_{m}) +L Q_{m+1}=&\\\frac{1}{2}\mq(\tau L)f(Q_{m})+\frac{1}{2}(\mq(\tau L)+\tau L)&f(e^{-\tau L}Q_{m})+\frac{1}{2}\tau \mq(\tau L)\frac{\partial f}{\partial Q}(Q_{m}):f(Q_{m}), \label{lri2a_2}\\
    \frac{1}{\tau} \mq(\tau L)(Q_{m+1}-Q_{m}) +L Q_{m+1}& =\\\frac{1}{2}\mq(\tau L)f(Q_{m})+\frac{1}{2}(\mq(\tau L)+\tau L)f(e^{-\tau L}Q_{m})&+\frac{1}{2}\tau (\mq(\tau L)+\tau L)\frac{\partial f}{\partial Q}(e^{-\tau L}Q_{m}):f(e^{-\tau L}Q_{m}). \label{lri2b_2}
\end{align}
Utilizing the fact that  $L$ is  positive semidefinite,   we define the following inner products and the corresponding energy norms for any $\mathbf{U} \in \mathcal{X}\cap \mathcal{Z}$ and  $\tau > 0$:
\begin{align}
 \ilQ{ \mathbf{U} }
 :=   (L\mathbf{U},\mathbf{U})_{\mathcal{X}} ,\quad  \ilQa{ \mathbf{U} }
 :=   (L_1\mathbf{U},\mathbf{U})_{\mathcal{X}} ,\quad
    \iaQ{ \mathbf{U} }=  (\psi( \tau L) \mathbf{U}, \mathbf{U})_{\mathcal{X}} ,  \\
\ibQ{ \mathbf{U} }= {(\psi_1(\tau L) \mathbf{U}, \mathbf{U})_{\mathcal{X}}} = {(\psi(\tau L) \mathbf{U}, \mathbf{U})_{\mathcal{X}}} +  {(\frac{\tau}{2} L \mathbf{U}, \mathbf{U})_{\mathcal{X}}} = \iaQ{ \mathbf{U} } + \frac{\tau}{2}\ilQ{ \mathbf{U} }.
\end{align}
where the  function $\psi_1$ is defined as
\begin{align}\label{eq_Q_functions}
    \psi_1(z) =\psi(z) + \frac{z}{2} = \frac{z}{e^z - 1}+ \frac{z}{2}= \frac{z}{2}\text{coth}\left( \frac{z}{2} \right)
    \geq 1, \quad   z\geq 0.
\end{align}
With the above definitions, we can define the discrete energy functional $\mathcal{E}$, $\mathcal{E}_i[Q]$  as follows:
\begin{align}
    	\mathcal{E}[Q]= \frac{1}{2}\ilQ{Q}+E[Q],\quad
        \mathcal{E}_1[Q]= \frac{1}{2}\ilQa{Q}+E[Q]. \label{discrete_energy}
\end{align}
where
\begin{align*}
    E[Q] &= \frac{\alpha}{2} \| Q\|_{\mathcal{X}}^2-\frac{\beta}{3} (Q,Q^2)_{\mathcal{X}}+\frac{\gamma}{4} (\|Q\|_F^2,\|Q\|_F^2)_\mathcal{X}.
\end{align*}
\begin{lemma}\label{3.4}
	 For any $Q_1, Q_2 \in \mathcal{X}\cap \mathcal{Z}$ satisfying $\|Q_1\|_\mathcal{Z}\leq \eta, \|Q_2\|_\mathcal{Z} \leq \eta$, for any $\boldsymbol{x} \in \bar{\Omega}$, { when $\kappa\geq  |\alpha|+\beta \eta+3 \gamma \eta^2$,} we can derive
	 \begin{align*}
		 &E[Q_1] - E[Q_2 ]+  (Q_1 - Q_2, f(Q_2))_{\mathcal{X}}  \leq \kappa\|Q_1 - Q_2\|_{\mathcal{X}}^2.
	 \end{align*}
 \end{lemma}
\begin{proof}
For any $Q_1, Q_2 \in \mathcal{X}\cap \mathcal{Z}$ satisfying $\|Q_1\|_\mathcal{Z}\leq \eta, \|Q_2\|_\mathcal{Z} \leq \eta$, we can derive
\begin{align}
    \frac{\gamma}{4} (\|Q_1\|_F^2,\|Q_1\|_F^2)_\mathcal{X}-\frac{\gamma}{4} (\|Q_2\|_F^2,\|Q_2\|_F^2)_\mathcal{X}\le \frac{\gamma}{2} ((\|Q_1\|_{F}^2+\|Q_2\|_{F}^2)Q_1,Q_1-Q_2)_{\mathcal{X}}.
\end{align}
For $E$, we have
\begin{align}
    &\quad E[Q_1] - E[Q_2]\no\\
    =& \frac{\alpha}{2} (\|Q_1\|_{\mathcal{X}}^2-\|Q_2\|_{\mathcal{X}}^2)-\frac{\beta}{3} ((Q_1,Q_1^2)-(Q_2,Q_2^2))_{\mathcal{X}}+\frac{\gamma}{4} ((\|Q_1\|_F^2,\|Q_1\|_F^2)_\mathcal{X}-(\|Q_2\|_F^2,\|Q_2\|_F^2)_\mathcal{X})\no\\
    \leq & \alpha (Q_1,Q_1-Q_2)_{\mathcal{X}}-\frac{\beta}{3} (Q_1-Q_2, Q_1^2+(Q_1+Q_2)Q_2)_{\mathcal{X}}+ \frac{\gamma}{2} ((\|Q_1\|_{F}^2+\|Q_2\|_{F}^2)Q_1,Q_1-Q_2)_{\mathcal{X}}.
\end{align}

    	By the definition of   $f(Q)$, we can derive
\begin{align}
            &\quad (Q_1 - Q_2, f(Q_2))_{\mathcal{X}}\no\\ &= -\alpha (Q_1 - Q_2, Q_2)_{\mathcal{X}} + \beta (Q_1 - Q_2, Q_2^2 - \frac{1}{3}\text{tr}(Q_2^2)I)_{\mathcal{X}} - \gamma \text{tr}(Q_2^2)(Q_1 - Q_2, Q_2)_{\mathcal{X}}\no\\
		&= -\alpha (Q_1 - Q_2, Q_2)_{\mathcal{X}} + \beta (Q_1 - Q_2, Q_2^2)_{\mathcal{X}}  -\gamma  (Q_1 - Q_2, \|Q_2\|_{F}^2Q_2)_{\mathcal{X}}.\no
	\end{align}
	By adding the above inequalities together, and using $\kappa\geq |\alpha|+\beta \eta+3 \gamma \eta^2$ and $\|Q_1\|_{\mathcal{Z}},\|Q_2\|_{\mathcal{Z}}\le \eta$, we  obtain
        	\begin{align*}
		&\quad E[Q_1] - E[Q_2] +  (Q_1 - Q_2, f(Q_2))_{\mathcal{X}} \\
		&= \alpha (Q_1 - Q_2,Q_1 - Q_2)_{\mathcal{X}} -\frac{\beta}{3}(Q_1 - Q_2, Q_1^2 + Q_2^2+Q_1Q_2-3Q_2^2)_{\mathcal{X}}\\
		&\quad+\frac{\gamma}{2} (\|Q_1\|_{F}^2Q_1+\|Q_2\|_{F}^2Q_1-2\|Q_2\|_{F}^2Q_2,Q_1-Q_2)_{\mathcal{X}} \\
		&\leq  \left(\alpha +\frac{\beta}{3}(\|Q_1\|_{\mathcal{Z}} +2 \|Q_2\|_{\mathcal{Z}})+\frac{\gamma}{2}(2\|Q_2\|_{\mathcal{Z}}^2+\|Q_1\|_{\mathcal{Z}}^2+\|Q_1\|_{\mathcal{Z}}\|Q_2\|_{\mathcal{Z}}) \right)\|Q_1-Q_2\|_{\mathcal{X}}^2	\\
		&\leq (|\alpha|+\beta \eta+2\gamma \eta^2) \|Q_1-Q_2\|_{\mathcal{X}}^2.
		\end{align*}
\end{proof}
\begin{theorem}
    For any $Q_0 \in\mathcal{X}\cap \mathcal{Z}$ satisfying $\|Q_0\|^2_\mathcal{Z}\leq \eta^2$, when $\tau\leq\tau_0, m \geq1$, the solution $Q_m$ generated by the schemes \eqref{lri1a} satisfies
    \begin{align}
        \mathcal{E}_1[Q_{m}]  &\leq\mathcal{E}_1[Q_{m-1}], \quad m \geq1, \label{energy_stability}
    \end{align}
  and the solution $Q_m$ generated by the schemes \eqref{lri1b} satisfies
    \begin{align}
        \mathcal{E}_1[e^{-\tau L}Q_{m}]  &\leq\mathcal{E}_1[e^{-\tau L}Q_{m-1}], \quad m \geq1, \label{energy_stability2}
    \end{align}
    where $\mathcal{E}_1$ is defined in \eqref{discrete_energy}.
\end{theorem}
\begin{proof}
    We first consider the LRI1a scheme \eqref{lri1a_3}. Taking the inner product of both sides of \eqref{lri1a_3} with $Q_{m+1}-Q_{m}$, we have
    \begin{align}
        \quad \frac{1}{\tau}(Q_{m+1}-Q_{m}, Q_{m+1}-Q_{m})_{\mathcal{X}}+(L_1 Q_{m+1}, Q_{m+1}-Q_{m})_{\mathcal{X}} = (f(Q_{m}), Q_{m+1}-Q_{m})_{\mathcal{X}}.
    \end{align}
    By Lemma \ref{3.4}, we can derive
    \begin{align}
        &\quad \mathcal{E}_1[Q_{m+1}] - \mathcal{E}_1[Q_{m}] \no\\
        &= \frac{1}{2}\ilQa{Q_{m+1}} - \frac{1}{2}\ilQa{Q_{m}}+E[Q_{m+1}] - E[Q_{m}] \no\\
        &\leq (L_1 Q_{m+1}, Q_{m+1}-Q_{m})_{\mathcal{X}}-\frac{1}{2}(L_1(Q_{m+1}-Q_{m}), Q_{m+1}-Q_{m})_{\mathcal{X}} +E[Q_{m+1}] - E[Q_{m}]\\
        &\leq (f(Q_{m}), Q_{m+1}-Q_{m})_{\mathcal{X}} +E[Q_{m+1}] - E[Q_{m}] -\frac{1}{\tau}\norm{Q_{m+1}-Q_{m}}_{\mathcal{X}}^2\\
        &\le \kappa \|Q_{m+1}-Q_{m}\|_{\mathcal{X}}^2 -\frac{1}{\tau}\norm{Q_{m+1}-Q_{m}}_{\mathcal{X}}^2\leq 0,
    \end{align}
    where we have used the fact that $\tau\leq \tau_0\leq \frac{1}{\kappa}$.

We demonstrate that the LRI1b scheme can be exactly transformed into the LRI1a scheme via a straightforward variable substitution. Recall the LRI1b scheme \eqref{lri1b}:
\begin{equation}\label{eq:lri1b}
    Q_{m+1} = e^{c\tau \Delta} Q_m + \tau f(e^{c\tau \Delta} Q_m).
\end{equation}
Introducing  $P_m = e^{c\tau \Delta} Q_m$ and applying  $e^{c\tau \Delta}$ to both sides of \eqref{eq:lri1b}, we obtain:
\begin{equation}\label{eq:sub_step2}
    e^{c\tau \Delta} Q_{m+1} = e^{c\tau \Delta} \left( P_m + \tau f(P_m) \right).
\end{equation}
By definition, the left-hand side is precisely $P_{m+1}$, leading to
\begin{equation}\label{eq:lri1a_derived}
    P_{m+1} = e^{c\tau \Delta} P_m + \tau e^{c\tau \Delta} f(P_m),
\end{equation}
which is exactly the LRI1a scheme. This illustrates that LRI1a and LRI1b fundamentally generate an identical alternating sequence of physical states, differing only in the sequence observation point (i.e., whether diffusion or reaction is performed first).

Consequently, as the LRI1a scheme satisfies the strictly monotonic modified energy dissipation law $\mathcal{E}_1(P_{m+1}) \le \mathcal{E}_1(P_m)$, the LRI1b scheme inherently preserves this stability.
\end{proof}
\begin{lemma}\label{lem3_4}
    For any $Q_0 \in\mathcal{X}\cap \mathcal{Z}$ satisfying $\|Q_0\|^2_\mathcal{Z}\leq \eta^2$, when $\tau\leq\tau_0, m \geq1$, the solution $Q_m$ generated by the schemes \eqref{lri1b} and \eqref{lri2b} satisfies
    \begin{align}
      \frac{1}{2}\|Q_M\|_{L}^2 + \sum_{m=0}^{M-1} \frac{1}{4\tau}\|Q_{m+1}-Q_m\|_{\mathcal{X}}^2\leq C(\norm{Q_0}_{L},\eta,T). \label{energy_stability3}
    \end{align}
\end{lemma}
\begin{proof}
Here we analyze the stability of \eqref{lri2}, and the stability of \eqref{lri1b} and \eqref{lri2b} can be derived similarly.

For brevity, we denote $\tilde{Q}_m = e^{-\tau L} Q_m$.
Taking the $L^2$ inner product of the scheme with the time difference test function $Q_{m+1} - Q_m$, we have:
\begin{equation}\label{eq:inner_prod_2nd}
    \begin{aligned}
        &\left(\mq(\tau L) \frac{Q_{m+1}-Q_{m}}{\tau}, Q_{m+1}-Q_{m} \right)_{\mathcal{X}} + (L Q_{m+1}, Q_{m+1}-Q_m)_{\mathcal{X}} \\
        =& \left( \frac{1}{2}\mq(\tau L)(f(Q_m) +f(\tilde{Q}_m)), Q_{m+1}-Q_m \right)_{\mathcal{X}} + \left( \frac{1}{2}\tau L f(\tilde{Q}_m), Q_{m+1}-Q_m \right)_{\mathcal{X}} \\
        &+ \left( \frac{1}{2}\tau \mq(\tau L)\frac{\partial f}{\partial Q}(Q_m):f(Q_m), Q_{m+1}-Q_m \right)_{\mathcal{X}}=:I_1+I_2+I_3.
    \end{aligned}
\end{equation}
Using the standard identity for the self-adjoint operator $L$, the LHS of \eqref{eq:inner_prod_2nd} can be expanded as follows:
\begin{align}\label{eq:lhs_expand}
    \text{LHS} &= \frac{1}{\tau}\|Q_{m+1}-Q_m\|_{\psi}^2 + \frac{1}{2}\|Q_{m+1}\|_{L}^2 - \frac{1}{2}\|Q_m\|_{L}^2 + \frac{1}{2}\|Q_{m+1}-Q_m\|_{L}^2,\\
    &= \frac{1}{2\tau}\|Q_{m+1}-Q_m\|_{\psi}^2 + \frac{1}{4}\|Q_{m+1}-Q_m\|_{L}^2 + \frac{1}{2}\|Q_{m+1}\|_{L}^2 - \frac{1}{2}\|Q_m\|_{L}^2\\
    &+\frac{1}{2\tau}\|Q_{m+1}-Q_m\|_{\psi}^2 + \frac{1}{4}\|Q_{m+1}-Q_m\|_{L}^2\\
    &= \frac{1}{2\tau}\|Q_{m+1}-Q_m\|_{\psi_L^1}^2  + \frac{1}{2}\|Q_{m+1}\|_{L}^2 - \frac{1}{2}\|Q_m\|_{L}^2+\frac{1}{2\tau}\|Q_{m+1}-Q_m\|_{\psi}^2 + \frac{1}{4}\|Q_{m+1}-Q_m\|_{L}^2
\end{align}
We estimate the three inner products on the RHS separately using the Cauchy-Schwarz and Young's inequalities. Assume the time step $\tau \le 1$.
For the standard nonlinear terms:
\begin{equation}\label{eq:rhs1}
    I_1 \le \frac{1}{8\tau}\|Q_{m+1}-Q_m\|_{\mathcal{X}}^2 + C\tau \left( \|f(Q_m)\|_{\mathcal{X}}^2 + \|f(\tilde{Q}_m)\|_{\mathcal{X}}^2 \right).
\end{equation}

For the stabilization term, we distribute the operator $L^{1/2}$:
\begin{equation}\label{eq:rhs2}
    \begin{aligned}
        I_2 &= \frac{\tau}{2} \left( L^{1/2} f(\tilde{Q}_m), L^{1/2}(Q_{m+1}-Q_m) \right) \\
        &\le \frac{1}{4}\|Q_{m+1}-Q_m\|_{L}^2 + \frac{\tau^2}{4} \|f(\tilde{Q}_m)\|_{L}^2 .
    \end{aligned}
\end{equation}

For the Jacobian-related term:
\begin{equation}\label{eq:rhs3}
    \begin{aligned}
        I_3 &\le \frac{1}{8\tau}\|Q_{m+1}-Q_m\|_{\mathcal{X}}^2 + 2\tau \left\| \frac{1}{2}\tau \frac{\partial f}{\partial Q}(Q_m):f(Q_m) \right\|_{\mathcal{X}}^2 \\
        &\le \frac{1}{8\tau}\|Q_{m+1}-Q_m\|_{\mathcal{X}}^2 + C\tau^3 \left\| \frac{\partial f}{\partial Q}(Q_m):f(Q_m) \right\|_{\mathcal{X}}^2 \\
        &\le \frac{1}{8\tau}\|Q_{m+1}-Q_m\|_{\mathcal{X}}^2 + C\tau \left\| \frac{\partial f}{\partial Q}(Q_m):f(Q_m) \right\|_{\mathcal{X}}^2.
    \end{aligned}
\end{equation}
Substituting \eqref{eq:rhs1}, \eqref{eq:rhs2}, and \eqref{eq:rhs3} into \eqref{eq:inner_prod_2nd}, and using the expansion of the LHS in \eqref{eq:lhs_expand}, we arrive at the following inequality:
\begin{equation}
    \begin{aligned}
        &\frac{1}{4\tau}\|Q_{m+1}-Q_m\|_{\psi_L^1}^2  + \frac{1}{2}\|Q_{m+1}\|_{L}^2 - \frac{1}{2}\|Q_m\|_{L}^2 \\
        &\le C\tau \left( \|f(Q_m)\|_{\mathcal{X}}^2 + \|f(\tilde{Q}_m)\|_{L}^2 + \left\| \frac{\partial f}{\partial Q}(Q_m):f(Q_m) \right\|_{\mathcal{X}}^2 \right).
    \end{aligned}
\end{equation}
Summing from $m=0$ to $M-1$ and utilizing the telescoping property of $\|Q_m\|_{L}^2$:
\begin{equation}\label{eq:telescoping}
    \begin{aligned}
        &\frac{1}{2}\|Q_M\|_{L}^2 + \sum_{m=0}^{M-1} \frac{1}{4\tau}\|Q_{m+1}-Q_m\|_{\mathcal{X}}^2 \\
        &\le \frac{1}{2}\|Q_0\|_{L}^2 + C\tau \sum_{m=0}^{M-1} \left( \|f(Q_m)\|_{\mathcal{X}}^2 + \|f(\tilde{Q}_m)\|_{L}^2 + \left\| \frac{\partial f}{\partial Q}(Q_m):f(Q_m) \right\|_{\mathcal{X}}^2 \right).
    \end{aligned}
\end{equation}

Applying  Theorem \ref{mbp}, all nonlinear terms on the RHS are strictly controlled by the $L$ norm of $Q_m$:
\begin{equation}
    \|f(Q_m)\|_{\mathcal{X}}^2 + \|f(\tilde{Q}_m)\|_{L}^2 + \left\| \frac{\partial f}{\partial Q}(Q_m):f(Q_m) \right\|_{\mathcal{X}}^2 \le C(\eta) \|Q_m\|_{L}^2.
\end{equation}
Dropping the positive time-difference summation on the LHS of \eqref{eq:telescoping} for the inequality, we obtain the standard canonical form for the discrete Gronwall's lemma:
\begin{equation}\label{eq:gronwall_ready}
    \|Q_M\|_{L}^2 \le \|Q_0\|_{L}^2 + 2C(\eta) \tau \sum_{m=0}^{M-1} \|Q_m\|_{L}^2.
\end{equation}
Applying the discrete Gronwall's inequality yields:
\begin{equation}\label{eq:gronwall_final}
    \|Q_M\|_{L}^2 \le \|Q_0\|_{L}^2 \exp\left( 2C(\eta) \tau M \right)= \|Q_0\|_{L}^2\exp(2C(\eta) T).
\end{equation}
Substituting the above bound back into \eqref{eq:telescoping}, we can derive the desired stability result.
\end{proof}
\begin{theorem}\label{thm:energy_stability}
    For any $Q_0 \in H^{1}\cap \mathcal{Z}$ satisfying $\|Q_0\|^2_\mathcal{Z}\leq \eta^2$,  $Q_m$  is generated by schemes \eqref{lri2}-\eqref{lri2b} when $\tau\leq\tau_0$.  We have the following uniform energy stability bound:
    \begin{align}
        \mathcal{E}[Q_m] \leq C(\norm{Q_0}_{L}, \eta, T, |\Omega|), \quad \forall m \ge 0 \text{ such that } m\tau \leq T,
    \end{align}
    where $C>0$ is a uniform constant independent of the time step $\tau$ and the spatial grid size.
\end{theorem}

\begin{proof}
    Recall that the total free energy $\mathcal{E}[Q_{m+1}]$ can be decoupled into the elastic gradient energy and the nonlinear bulk energy $E[Q_{m+1}] = \int_{\Omega} F(Q_{m+1}) dx$.
        Theorem \ref{mbp} guarantees that the bulk energy $E[Q_{m+1}]$ is uniformly bounded by a constant depending on $\eta$ and the domain size $|\Omega|$:
    \begin{align}\label{eq:bulk_bound}
        E[Q_{m+1}] \le \int_{\Omega} \left| F(Q_{m+1}) \right| dx \leq C_F(\eta)|\Omega|.
    \end{align}
    Combining the $H^1$ bound from Lemma \ref{lem3_4} and the bulk energy bound \eqref{eq:bulk_bound} derived from Theorem \ref{mbp}, we can complete the derivation for the total energy:
    \begin{align}
        \mathcal{E}[Q_{m+1}] &\leq \frac{1}{2}\ilQ{Q_{m+1}} + E[Q_{m+1}]  \no\\
        &\leq  C(\norm{Q_0}_{L},\eta,T) + C_F(\eta)|\Omega|.
    \end{align}
    Letting $C(\norm{Q_0}_{L}, \eta, T, |\Omega|) := C(\norm{Q_0}_{L},\eta,T) + C_F(\eta)|\Omega|$ yields the desired uniform upper bound for the total energy. This completes the proof.
\end{proof}
 \section{Temporal error estimate}\label{section4}~\\
In this section, we provide a rigorous temporal error analysis for schemes \eqref{lri1a}--\eqref{lri2b}.
The core of our approach is to leverage the first and second derivatives of $Y(\xi,Q)$ to effectively bound the temporal error arising from the nonlinearities.
Here , we define the auxiliary function $Y(\xi,Q)$ and analyze its properties.
 \begin{lemma}\label{th3.2}
 	For each $0\leq \xi\leq \tau$,  the function $Y(\xi,Q):\mathcal{R}\times \mathcal{X}\to \mathcal{X}$ is defined as
 	\begin{align}
 		Y(\xi,Q)=e^{c(\tau-\xi) \Delta} f(e^{c\xi \Delta}Q). \label{3.6}
 	\end{align}
	 $Y'(\xi,Q)$ and $Y''(\xi,Q)$ are denoted as the first and second derivatives of $Y(\xi,Q)$ with respect to $\xi$.
Then for any given $\xi$, we have:
\begin{enumerate}
    \item[(i)] For $Q \in H^1(\Omega)$, when $f$ is second-order differentiable,  $Y'(\xi,Q)$ is bounded in the $\mathcal{X}$-norm.

    \item[(ii)] For $Q \in H^2(\Omega)$, when $f$ is forth-order differentiable,  $Y''(\xi,Q)$ is bounded in the $\mathcal{X}$-norm.
\end{enumerate}
 	 \end{lemma}
 	\begin{proof}
 		Letting $M=e^{c\xi \Delta}Q$, for $Y'(\xi,Q)_{ij}$, we have
 			\begin{align}
 			Y'(\xi,Q)_{ij}&= e^{c(\tau-\xi) \Delta}\left(-c \Delta f_{ij}(M)+\frac{\partial f_{ij}}{\partial Q_{lm}}(M) c\Delta  M_{lm}\right)\no\\&=
 			e^{c(\tau-\xi)\Delta}\left(-c \nabla   \cdot \nabla f_{ij}(M)+\frac{\partial f_{ij}}{\partial Q_{lm}}(M) c\Delta  M_{lm}\right)
 			\no\\&=
 			e^{c(\tau-\xi)\Delta}\left(-c \nabla   \cdot\left(\frac{\partial f_{ij}}{\partial Q_{lm}}(M)  \frac{\partial M_{lm}}{\partial x_k}\right) +\frac{\partial f_{ij}}{\partial Q_{lm}}(M) c\Delta  M_{lm}\right)\no\\&=
 			e^{c(\tau-\xi)\Delta}\left(-c \left(\frac{\partial(\frac{\partial f_{ij}}{\partial Q_{lm}}(M))}{\partial x_k}  \frac{\partial M_{lm}}{\partial x_k}+\frac{\partial f_{ij}}{\partial Q_{lm}}(M)  \frac{\partial^2 M_{lm}}{\partial x_k\partial x_k}\right) +\frac{\partial f_{ij}}{\partial Q_{lm}}(M) c\Delta  M_{lm}\right)\no\\&=
 			e^{c(\tau-\xi)\Delta}\left(-c \frac{\partial^2 f_{ij}}{\partial Q_{lm}\partial Q_{st}}(M) \frac{\partial M_{st}}{\partial x_k} \frac{\partial M_{lm}}{\partial x_k}\right).\no
 		\end{align}
		 Taking the $\mathcal{X}$ norm on both sides of the above equation and using $Q \in H^1$, we can obtain $\Norm{Y'(\xi,Q)}{\mathcal{X}}$ is bounded.

		Letting $f'',f''',f''''$  denote $\frac{\partial^2 f_{ij}}{\partial Q_{lm}\partial Q_{st}}(M),\frac{\partial^3 f_{ij}}{\partial Q_{lm}\partial Q_{st}\partial Q_{up}}(M),\frac{\partial^4 f_{ij}}{\partial Q_{lm}\partial Q_{st}\partial Q_{up}\partial Q_{ow}}(M)$,
		for the term $Y''(\xi,Q)_{ij}$, we have
		\begin{align*}
			Y''(\xi,Q)_{ij}&=
			e^{c(\tau-\xi)\Delta}\left(c^2 \Delta (f''\frac{\partial M_{st}}{\partial x_k} \frac{\partial M_{lm}}{\partial x_k})-c f''c\Delta\frac{\partial M_{st}}{\partial x_k} \frac{\partial M_{lm}}{\partial x_k}-c f''\frac{\partial M_{st}}{\partial x_k} c\Delta\frac{\partial M_{lm}}{\partial x_k}\right)\\
		\end{align*}
For the first term of  $Y''(\xi,Q)_{ij}$, we have
\begin{align*}
	&c^2 \Delta (f''\frac{\partial M_{st}}{\partial x_k} \frac{\partial M_{lm}}{\partial x_k}) \\&=c^2 \nabla \cdot (f'''\frac{\partial M_{up}}{\partial x_g}\frac{\partial M_{st}}{\partial x_k} \frac{\partial M_{lm}}{\partial x_k}+f''\frac{\partial^2 M_{st}}{\partial x_k\partial x_g} \frac{\partial M_{lm}}{\partial x_k}+f''\frac{\partial M_{st}}{\partial x_k} \frac{\partial^2 M_{lm}}{\partial x_k\partial x_g})\\
	&=c^2(f''''\frac{\partial M_{ow}}{\partial x_g}\frac{\partial M_{up}}{\partial x_g}\frac{\partial M_{st}}{\partial x_k} \frac{\partial M_{lm}}{\partial x_k}\\&
	+f'''\frac{\partial M_{up}}{\partial x_g}\frac{\partial M_{st}}{\partial x_k} \frac{\partial^2 M_{lm}}{\partial x_k\partial x_g}+f'''\frac{\partial M_{up}}{\partial x_g}\frac{\partial^2 M_{st}}{\partial x_k\partial x_g} \frac{\partial M_{lm}}{\partial x_k}+f'''\frac{\partial^2 M_{up}}{\partial x_g\partial x_g}\frac{\partial M_{st}}{\partial x_k} \frac{\partial M_{lm}}{\partial x_k}\\
	&+f'''\frac{\partial M_{up}}{\partial x_g}\frac{\partial^2 M_{st}}{\partial x_k\partial x_g} \frac{\partial M_{lm}}{\partial x_k}+f''\frac{\partial^3 M_{st}}{\partial x_k\partial x_g\partial x_g} \frac{\partial M_{lm}}{\partial x_k}+f''\frac{\partial^2 M_{st}}{\partial x_k\partial x_g} \frac{\partial^2 M_{lm}}{\partial x_k\partial x_g}\\
	&+f'''\frac{\partial M_{up}}{\partial x_g}\frac{\partial M_{st}}{\partial x_k} \frac{\partial^2 M_{lm}}{\partial x_k\partial x_g}+f''\frac{\partial^2 M_{st}}{\partial x_k\partial x_g} \frac{\partial^2 M_{lm}}{\partial x_k\partial x_g}+f''\frac{\partial M_{st}}{\partial x_k} \frac{\partial^3 M_{lm}}{\partial x_k\partial x_g\partial x_g}).
\end{align*}
Then we have
\begin{align*}
	Y''(\xi,Q)_{ij}&=	c^2 e^{c(\tau-\xi)\Delta}(f''''\frac{\partial M_{ow}}{\partial x_g}\frac{\partial M_{up}}{\partial x_g}\frac{\partial M_{st}}{\partial x_k} \frac{\partial M_{lm}}{\partial x_k}+f''\frac{\partial^2 M_{st}}{\partial x_k\partial x_g} \frac{\partial^2 M_{lm}}{\partial x_k\partial x_g}\\&
	+f'''\frac{\partial M_{up}}{\partial x_g}\frac{\partial M_{st}}{\partial x_k} \frac{\partial^2 M_{lm}}{\partial x_k\partial x_g}+f'''\frac{\partial M_{up}}{\partial x_g}\frac{\partial^2 M_{st}}{\partial x_k\partial x_g} \frac{\partial M_{lm}}{\partial x_k}+f'''\frac{\partial^2 M_{up}}{\partial x_g\partial x_g}\frac{\partial M_{st}}{\partial x_k} \frac{\partial M_{lm}}{\partial x_k}\\
	&+f'''\frac{\partial M_{up}}{\partial x_g}\frac{\partial^2 M_{st}}{\partial x_k\partial x_g} \frac{\partial M_{lm}}{\partial x_k}+f''\frac{\partial^2 M_{st}}{\partial x_k\partial x_g} \frac{\partial^2 M_{lm}}{\partial x_k\partial x_g}
	+f'''\frac{\partial M_{up}}{\partial x_g}\frac{\partial M_{st}}{\partial x_k} \frac{\partial^2 M_{lm}}{\partial x_k\partial x_g}).
\end{align*}
 	Taking the $\mathcal{X}$ norm on both sides of the above equation and using $Q \in H^2$, we have $\Norm{Y''(\xi,Q)}{\mathcal{X}}$ is bounded.
 	\end{proof}
 \begin{theorem}\label{theorem_lri_1a}
	For  a fixed terminal time $T$,
 	assume that $Q\in C([0,T];H^{1}\cap \mathcal{Z} )$ is the exact solution of \eqref{1.6a}.  $Q_m$ is the numerical solution generated by the LRI1a scheme \eqref{lri1a} and we set  $e_m=Q(t_m)-Q_m$. Then, we have
 	\begin{align}
 		\|Q_m-Q(t_m)\|_\mathcal{X}\leq \frac{C_{5}(e^{CT}-1)}{C}\tau,
 	\end{align}
 	for $\tau\leq \tau_0, m\geq 1$, where $C_{5}$  is a constant that is independent of $\tau$.
 \end{theorem}

 \begin{proof}
 	For the LRI1a scheme \eqref{lri1a}, we have
\begin{align}
	Q(t_{m+1})=e^{c\tau \Delta} Q(t_m)+\tau e^{c\tau \Delta}f(Q(t_m))+R_1(t_m),\label{3.5}
\end{align}
where $R_1(t_m)$ is the corresponding truncation error.
Subtracting \eqref{3.5} from \eqref{lri1a}, we derive
\begin{align}
	e_{m+1}=e^{c\tau \Delta} e_{m}+\tau e^{c\tau \Delta}(f(Q(t_m)-f(Q_m)))+R_1(t_m),
\end{align}
Taking the  norm  $\|\cdot\|_\mathcal{X}$ on both sides of the above equation and using Lemma \ref{1.17}, we have
\begin{align}
	\Norm{e_{m+1}}{\mathcal{X}}&=\Norm{ e_{m}}{\mathcal{X}}+\tau \Norm{(f(Q(t_m)-f(Q_m)))}{\mathcal{X}}+\Norm{R_1(t_m)}{\mathcal{X}},\no\\
&\leq \Norm{ e_{m}}{\mathcal{X}}+\tau C \Norm{e_{m}}{\mathcal{X}}+\Norm{R_1(t_m)}{\mathcal{X}}\label{em}
\end{align}

By comparing \eqref{1.9} and \eqref{3.5}, for $R_1(t_m)$, we can derive
	\begin{align}
		R_1(t_m)&=\int_{0}^{\tau}e^{c(\tau-\xi) \Delta}f(Q(t_m+\xi))d\xi-\tau e^{c\tau \Delta}f(Q(t_m))\no\\
		&=\int_{0}^{\tau}e^{c(\tau-\xi) \Delta}f(Q(t_m+\xi))- e^{c\tau \Delta}f(Q(t_m))d\xi\no\\
		&=\int_{0}^{\tau}e^{c(\tau-\xi) \Delta}f(Q(t_m+\xi))-e^{c(\tau-\xi) \Delta}f(e^{c\xi \Delta}Q(t_m))\no\\
		&\quad+e^{c(\tau-\xi) \Delta}f(e^{c\xi \Delta}Q(t_m))- e^{c\tau \Delta}f(Q(t_m))d\xi,\label{4.7}
	\end{align}
Taking the norm  $\|\cdot\|_\mathcal{X}$ on both sides of \eqref{4.7} and using  Lemma \ref{th3.2}, we have
\begin{align}
	&\Norm{R_1(t_m)}{\mathcal{X}} \\&\leq \int_{0}^{\tau}\Norm{f(Q(t_m+\xi))-f(e^{c\xi \Delta}Q(t_m))}{\mathcal{X}}+\Norm{e^{c(\tau-\xi) \Delta}f(e^{c\xi \Delta}Q(t_m))- e^{c\tau \Delta}f(Q(t_m))}{\mathcal{X}}d\xi\no\\
	&=\int_{0}^{\tau} C\Norm{ \int_{0}^{\xi} e^{c(\xi-\sigma)\Delta}f(Q(t_m+\sigma)) d\sigma}{\mathcal{X}}+\Norm{\int_{0}^{\xi} Y'(\xi,Q)d\xi}{\mathcal{X}} d \xi\no\\&\leq \int_{0}^{\tau} CC_f a\xi+  \Norm{  Y'(\xi,Q)}{\mathcal{X}}\xi d \xi\no\\
	&\leq\frac{1}{2}(CC_f a+  \Norm{  Y'(\xi,Q)}{\mathcal{X}})\tau^2
	:\leq C_{5} \tau^2\label{3.8}.
\end{align}
Substituting \eqref{3.8} into \eqref{em}, we derive
\begin{align}
	\Norm{e_{m+1}}{\mathcal{X}}&\leq \Norm{ e_{m}}{\mathcal{X}}+\tau C\Norm{ e_{m}}{\mathcal{X}}+C_{5}\tau^2, \label{3.9}
\end{align}
which implies that
\begin{align}
	\Norm{e_{m}}{\mathcal{X}}+\frac{C_{5}\tau}{C}&\leq (1+C\tau)(\Norm{ e_{m-1}}{\mathcal{X}}+\frac{C_{5}\tau}{C})\no\\&\leq(1+C\tau)^{m}(\Norm{ e_{0}}{\mathcal{X}}+\frac{C_{5}\tau}{C}).
\end{align}
Since $\Norm{e_0}{\mathcal{X}}=0$, we have
\begin{align}
	\Norm{e_{m}}{\mathcal{X}}\leq((1+C\tau)^{m}-1)(\frac{C_{5}\tau}{C})
	\leq \frac{C_{5}(e^{CT}-1)}{C}\tau.\label{4.9}
\end{align}
 \end{proof}
 \begin{theorem}\label{theorem_lri_1b}
	For  a fixed terminal time $T$,
 	assume  that $Q\in C([0,T];H^{1}\cap \mathcal{Z})$ is the exact solution of \eqref{1.6a} and that $Q_m$ is the numerical solution generated by the LRI1b scheme \eqref{lri1b}. Then,   we have
 	\begin{align}
 		\|Q_m-Q(t_m)\|_\mathcal{X}\leq \frac{C_{5}(e^{CT}-1)}{C}\tau,\label{4.10}
 	\end{align}
	for $\tau\leq \tau_0, m\geq 1$,
 	where $C_{5}$ is a constant that is independent of $\tau$.
 \end{theorem}

 \begin{proof}
 	For the LRI1b scheme \eqref{lri1b}, we have
 	\begin{align}
 		Q(t_{m+1})=e^{c\tau \Delta} Q(t_m)+\tau f(e^{c\tau \Delta}Q(t_m))+R_2(t_m),\label{3.5a}
 	\end{align}
 	where $R_2(t_m)$ is the corresponding truncation error.
	Subtracting \eqref{lri1b} from \eqref{3.5a}, we derive
 	\begin{align}
 		Q(t_{m+1})-Q_m&= e^{c\tau \Delta} (Q(t_m)-Q_m)+\tau (f(e^{c\tau \Delta}Q(t_m))-f(Q_m))+R_2(t_m).\no\\
		 e_{m+1}&=e^{c\tau \Delta} e_{m}+\tau (f(e^{c\tau \Delta}Q(t_m)-f(e^{c\tau \Delta}Q_m)))+R_2(t_m).\label{3.5b}
 	\end{align}
	Taking the $\mathcal{X}$ norm  on both sides of \eqref{3.5b} and using Lemma \ref{1.17}, we derive
 	\begin{align}
 		\Norm{e_{m+1}}{\mathcal{X}}=\Norm{ e_{m}}{\mathcal{X}}+\tau C \Norm{Q(t_m)-Q_m}{\mathcal{X}}+\Norm{R_2(t_m)}{\mathcal{X}}.\no
 	\end{align}

 	On the other hand, by comparing \eqref{1.9} and \eqref{3.5a}, we derive
	\begin{alignat}{2}
		R_2(t_m)&=\int_{0}^{\tau}e^{c(\tau-\xi) \Delta}f(Q(t_m+\xi))d\xi-\tau f(e^{c\tau \Delta}Q(t_m))\no\\
			&=\int_{0}^{\tau}e^{c(\tau-\xi) \Delta}f(Q(t_m+\xi))- f(e^{c\tau \Delta}Q(t_m))d\xi\no\\
			&=\int_{0}^{\tau}e^{c(\tau-\xi) \Delta}f(Q(t_m+\xi))-e^{c(\tau-\xi) \Delta}f(e^{c\xi \Delta}Q(t_m))\no\\
			&\quad+e^{c(\tau-\xi) \Delta}f(e^{c\xi \Delta}Q(t_m))- f(e^{c\tau \Delta}Q(t_m))d\xi. \label{4-14}
	\end{alignat}
	Taking the $\mathcal{X}$ norm on both sides of \eqref{4-14} and using Lemma \ref{th3.2}, we have
	\begin{align*}
		&\Norm{R_2(t_m)}{\mathcal{X}}	\\&\leq \int_{0}^{\tau}\Norm{f(Q(t_m+\xi))-f(e^{c\xi \Delta}Q(t_m))}{\mathcal{X}}+\Norm{e^{c(\tau-\xi) \Delta}f(e^{c\xi \Delta}Q(t_m))- f(e^{c\tau \Delta}Q(t_m))d\xi}{\mathcal{X}}\no\\
		&\leq\int_{0}^{\tau} CC_f a\xi+\Norm{ \int_{\xi}^{\tau} Y'(\xi,Q) d\xi}{\mathcal{X}} d \xi\\&\leq \int_{0}^{\tau} CC_f a\xi+ \Norm{  Y'(\xi,Q)}{\mathcal{X}}(\tau -\xi) d \xi\no\\
		&\leq C_{5} \tau^2.\no
	\end{align*}
Similar to \eqref{3.9}-\eqref{4.9}, we obtain the conclusion \eqref{4.10}.
 \end{proof}

 \begin{theorem}\label{theorem_lri_2}
	For  a fixed terminal time $T$,
 	assume  that $Q\in C([0,T];H^{2}\cap \mathcal{Z})$ is the exact solution of \eqref{1.6a}. $Q_m$ represents the numerical solutions generated by the LRI2a scheme \eqref{lri2}, and $\overline{Q_m}$ is generated by the LRI2b scheme \eqref{lri2b}.
	Then, we have
 	\begin{align}
 		\|Q_m-Q(t_m)\|_\mathcal{X}&\leq \frac{C_{6}(e^{(C+\frac{1}{2}C_1\tau_0)T}-1)}{(C+\frac{1}{2}C_1\tau_0)}\tau^2,\no\\
		 \|\overline{Q_m}-Q(t_m)\|_\mathcal{X}&\leq \frac{C_{6}(e^{(C+\frac{1}{2}C_1\tau_0)T}-1)}{(C+\frac{1}{2}C_1\tau_0)}\tau^2,\no
 	\end{align}
	for $\tau\leq \tau_0$, $m\geq 1$ and $t_m\leq T$,
 	where $C_{6}$ is a constant that is independent of $\tau$..
 \end{theorem}

 \begin{proof}
 	For the LRI2a scheme \eqref{lri2}, we have
 	\begin{align}
 		Q(t_{m+1})&=e^{c\tau \Delta} Q(t_m)+\frac{1}{2}\tau e^{c\tau \Delta}f(Q(t_m))+\frac{1}{2}\tau f(e^{c\tau \Delta}Q(t_m))\no\\&+\frac{1}{2}\tau^2 e^{c\tau \Delta}\frac{\partial f}{\partial Q}(Q(t_m)):f(Q(t_m))+R_3(t_m),\label{3.5ab}
 	\end{align}
 	where $R_3(t_m)$ is the corresponding truncation error.
	Subtracting \eqref{lri2} from \eqref{3.5ab}, we have
 	\begin{align}
 		e_{m+1}&=e^{c\tau \Delta} e_{m}+\frac{1}{2}\tau (e^{c\tau \Delta}f(Q(t_m)-e^{c\tau \Delta}f(Q_m)))+\frac{1}{2}\tau (f(e^{c\tau \Delta}Q(t_m)-f(e^{c\tau \Delta}Q_m)))\no\\&+\frac{1}{2}\tau^2 e^{c\tau \Delta}(\frac{\partial f}{\partial Q}(Q(t_m)):f(Q(t_m))-\frac{\partial f}{\partial Q}(Q_m):f(Q_m))+R_3(t_m),
 	\end{align}
	Taking the norm $\|\cdot\|_\mathcal{X}$ on both sides of the above equation and using Lemma \ref{1.17}, we have
 	\begin{align}
 		\Norm{e_{m+1}}{\mathcal{X}}=\Norm{ e_{m}}{\mathcal{X}}+\tau C \Norm{Q(t_m)-Q_m}{\mathcal{X}}+\frac{1}{2}\tau^2 C_1 \Norm{Q(t_m)-Q_m}{\mathcal{X}}+\Norm{R_3(t_m)}{\mathcal{X}},\no\\
 		\Norm{e_{m+1}}{\mathcal{X}}\leq (1+C \tau+\frac{1}{2}C_1\tau^2 )\Norm{ e_{m}}{\mathcal{X}}+\Norm{R_3(t_m)}{\mathcal{X}}.\label{4-17}
 	\end{align}

 	On the other hand, by comparing \eqref{1.9} and \eqref{3.5ab}, we can derive
\begin{align}
    R_3(t_m)
    ={}& \int_{0}^{\tau} e^{c(\tau-\xi) \Delta} f(Q(t_m+\xi)) \,d\xi
    - \frac{\tau}{2} f(e^{c\tau \Delta}Q(t_m)) \notag \\
    & - \frac{\tau}{2} e^{c\tau \Delta}f(Q(t_m))
    - \frac{\tau^2}{2} e^{c\tau \Delta} \frac{\partial f}{\partial Q}(Q(t_m)) : f(Q(t_m)) \notag \\
    ={}& \int_{0}^{\tau} \Bigg[ e^{c(\tau-\xi) \Delta}f(Q(t_m+\xi))
    - \xi e^{c\tau \Delta} \frac{\partial f}{\partial Q}(Q(t_m)) : f(Q(t_m)) \notag \\
    & \qquad - \left(1-\frac{\xi}{\tau}\right) e^{c\tau \Delta}f(Q(t_m))
    - \frac{\xi}{\tau} f(e^{c\tau \Delta}Q(t_m)) \Bigg] \,d\xi \notag \\
    ={}& \int_{0}^{\tau} (L_1 + L_2 + L_3 + L_4) \,d\xi, \label{3.7ab}
\end{align}
 	where
    \begin{small}
       	\begin{align*}
  	&L_1= e^{c(\tau-\xi) \Delta}\left[ f(Q(t_m+\xi))-f(e^{c\xi \Delta}Q(t_m)+\xi f(e^{c\xi \Delta}Q(t_m)))\right],\\
 	&L_2=e^{c(\tau-\xi) \Delta}\left[f(e^{c\xi \Delta}Q(t_m)+\xi f(e^{c\xi \Delta}Q(t_m)))-f(e^{c\xi \Delta}Q(t_m))-\xi \frac{\partial f}{\partial Q}(e^{c\xi \Delta} Q(t_m)):f(e^{c\xi \Delta}Q(t_m))\right],\\
 	&L_3=\xi e^{c(\tau-\xi) \Delta}\frac{\partial f}{\partial Q}(e^{c\xi \Delta} Q(t_m)):f(e^{c\xi \Delta}Q(t_m))-\xi  e^{c \tau \Delta}\frac{\partial f}{\partial Q}( Q(t_m)):f(Q(t_m)),\\
 	&L_4=e^{c(\tau-\xi) \Delta}f(e^{c\xi \Delta}Q(t_m))-(1-\frac{\xi}{\tau})e^{c\tau \Delta}f(Q(t_m))-\frac{\xi}{\tau} f(e^{c\tau \Delta}Q(t_m)).
 	\end{align*}
    \end{small}

According to the definition of $Y(\xi,Q)$ in Lemma \ref{th3.2},  upon
taking the norm $\|\cdot\|_\mathcal{X}$  with respect  to  $L_1$, we obtain
		\begin{align}
			\Norm{L_1}{\mathcal{X}}&= \Norm{e^{c(\tau-\xi) \Delta}\left[ f(Q(t_m+\xi))-f(e^{c\xi \Delta}Q(t_m)+\xi f(e^{c\xi \Delta}Q(t_m)))\right]}{\mathcal{X}}\no\\
			&\leq C \Norm{ Q(t_m+\xi)-e^{c\xi \Delta}Q(t_m)-\xi f(e^{c\xi \Delta}Q(t_m))}{\mathcal{X}} \no\\
			&\leq C \Norm{ \int_{0}^{\xi} e^{c(\xi-\sigma)\Delta}f(Q(t_m+\sigma)) d\sigma-\xi f(e^{c\xi \Delta}Q(t_m))}{\mathcal{X}}\no\\
			&\leq C \Norm{ \int_{0}^{\xi} e^{c(\xi-\sigma)\Delta}f(Q(t_m+\sigma))-e^{c(\xi-\sigma)\Delta}f(e^{c\sigma \Delta}Q(t_m)) d\sigma}{\mathcal{X}}\no\\
			&+ C \Norm{ \int_{0}^{\xi} e^{c(\xi-\sigma)\Delta}f(e^{c\sigma \Delta}Q(t_m))-f(e^{c\xi \Delta}Q(t_m)) d\sigma}{\mathcal{X}}\no\\
			&\leq\frac{1}{2}C^2C_f \Norm{ Q(t_m)}{\mathcal{X}}\xi^2+  C \int_{0}^{\xi} \int_{\sigma}^{\xi} \Norm{ Y'(\varsigma)}{\mathcal{X}} d\varsigma d\sigma\no\\
			&\leq\frac{1}{2}C^2C_f \Norm{ Q(t_m)}{\mathcal{X}}\xi^2+\frac{1}{2} C \xi^2 \Norm{ Y'(\varsigma)}{\mathcal{X}}.\label{4-18}
		\end{align}

	From \eqref{lem2.2a}, by taking the norm $\|\cdot\|_\mathcal{X}$ with respect to   $L_2$, we derive
 	\begin{align}
 		\Norm{L_2}{\mathcal{X}}\leq \xi^2 \Norm{\frac{\partial^2 f}{\partial Q^2}(e^{c\xi \Delta} Q(t_m)):f(e^{c\xi \Delta}Q(t_m)):f(e^{c\xi \Delta}Q(t_m))}{\mathcal{X}}\leq\xi^2 C_fC_{\partial}\Norm{ Q(t_m)}{\mathcal{X}}^2.\label{4-20}
 	\end{align}

We	set $\widehat{Y}(\xi)=e^{c(\tau-\xi) \Delta}\frac{\partial f}{\partial Q}(e^{c\xi \Delta} Q(t_m)):f(e^{c\xi \Delta}Q(t_m))$.  Taking the norm $\|\cdot\|_\mathcal{X}$ with respect  to  $L_3$ and using Lemma \ref{th3.2} , we have
 	\begin{align}
 		\Norm{L_3}{\mathcal{X}}&\leq \xi \Norm{ e^{c(\tau-\xi) \Delta}\frac{\partial f}{\partial Q}(e^{c\xi \Delta} Q(t_m)):f(e^{c\xi \Delta}Q(t_m))-  e^{c \tau \Delta}\frac{\partial f}{\partial Q}( Q(t_m)):f(Q(t_m))}{\mathcal{X}}\no\\
 		&\leq \xi \int_{0}^{\xi} \Norm{ \widehat{Y}'(\sigma)}{\mathcal{X}}  d\sigma \leq   \Norm{ \widehat{Y}'(\sigma)}{\mathcal{X}}\xi^2.\label{4-21}
 	\end{align}

	 According to the definition of $Y(\xi,Q)$  in \eqref{3.6}, we can rewrite $L_4$ as follows:
 	\begin{align}
 		L_4=Y(\xi,Q)-(1-\frac{\xi}{\tau})Y(0)-\frac{\xi}{\tau} Y(\tau).\no
 	\end{align}
	 Letting $Y(0),Y(\tau)$ be Taylor expanded at $\xi$ and taking norm $\|\cdot\|_\mathcal{X}$ with respect  to $L_4$, we can obtain
 	\begin{align}
 		\Norm{L_4}{\mathcal{X}}\leq \frac{1}{2} \Norm{Y''(\varsigma)}{\mathcal{X}}(\xi\tau-\xi^2),\qquad \varsigma \in[0,\tau].\label{4-22}
 	\end{align}

	Taking the $\mathcal{X}$ norm on both sides of \eqref{3.7ab} and substituting   inequalities \eqref{4-18}-\eqref{4-22} into it, we can derive
 	\begin{align}
 		\Norm{R_3(t_m)}{\mathcal{X}}&\leq
 		\int_{0}^{\tau} \Norm{L_1}{\mathcal{X}}+\Norm{L_2}{\mathcal{X}}+\Norm{L_3}{\mathcal{X}}+\Norm{L_4}{\mathcal{X}}d\xi
 		\no\\&\leq \int_{0}^{\tau}\frac{1}{2}C^2C_f \Norm{ Q(t_m)}{\mathcal{X}}\xi^2+\frac{1}{2} C \xi^2 \Norm{ Y'(\varsigma)}{\mathcal{X}}+\xi^2 C_fC_{\partial}\Norm{ Q(t_m)}{\mathcal{X}}^2\no\\&~~+\Norm{ \widehat{Y}'(\sigma)}{\mathcal{X}}\xi^2+ \frac{1}{2} \Norm{Y''(\varsigma)}{\mathcal{X}}(\xi\tau-\xi^2)d\xi\no\\&:\leq C_{6}\tau^3.\label{3.8ab}
 	\end{align}
 	By substituting \eqref{3.8ab} into \eqref{4-17}, we have
 	\begin{align}
 		\Norm{e_{m}}{\mathcal{X}}&\leq (1+C \tau+\frac{1}{2}C_1\tau^2 )\Norm{ e_{m-1}}{\mathcal{X}}+C_{6}\tau^3\no \\&\leq(1+(C+\frac{1}{2}C_1\tau_0) \tau )\Norm{ e_{m-1}}{\mathcal{X}}+C_{6}\tau^3,\no\\
 		\Norm{e_{m}}{\mathcal{X}}+\frac{C_{6}\tau^2}{(C+\frac{1}{2}C_1\tau_0)}&\leq (1+(C+\frac{1}{2}C_1\tau_0)\tau)(\Norm{ e_{m-1}}{\mathcal{X}}+\frac{C_{6}\tau^2}{(C+\frac{1}{2}C_1\tau_0)})\no\\&\leq(1+(C+\frac{1}{2}C_1\tau_0)\tau)^{m}(\Norm{ e_{0}}{\mathcal{X}}+\frac{C_{6}\tau^2}{(C+\frac{1}{2}C_1\tau_0)}).\label{label}
 	\end{align}
 	Substituting $\Norm{e_0}{\mathcal{X}}=0$ in \eqref{label}, we can drive
 	\begin{align*}
 		\Norm{e_{m}}{\mathcal{X}}&\leq((1+(C+\frac{1}{2}C_1\tau_0)\tau)^{m}-1)(\Norm{ e_{0}}{\mathcal{X}}+\frac{C_{6}\tau^2}{(C+\frac{1}{2}C_1\tau_0)})\\
 		&\leq \frac{C_{6}(e^{(C+\frac{1}{2}C_1\tau_0)T}-1)}{(C+\frac{1}{2}C_1\tau_0)}\tau^2.
 	\end{align*}

	 For the LRI2b scheme \eqref{lri2b},  the difference in the proof is
	 \begin{align*}
		\Norm{L_3}{\mathcal{X}}&\leq \xi \Norm{ e^{c(\tau-\xi) \Delta}\frac{\partial f}{\partial Q}(e^{c\xi \Delta} Q(t_m)):f(e^{c\xi \Delta}Q(t_m))-  \frac{\partial f}{\partial Q}( e^{c \tau \Delta}Q(t_m)):f(e^{c \tau \Delta}Q(t_m))}{\mathcal{X}}\\
		&\leq \xi \int_{\xi}^{\tau} \Norm{ \widehat{Y}'(\sigma)}{\mathcal{X}}  d\sigma \leq   \Norm{ \widehat{Y}'(\sigma)}{\mathcal{X}}\xi(\tau-\xi).
	\end{align*}
The remainder of the proof is similar to the above processes, and we omit the details.
 \end{proof}
 \begin{remark}
Notably, the proof presented above utilizes only the property that the boundary terms vanish during integration by parts. Thus, the same conclusion can also be obtained when  homogeneous Dirichlet boundary conditions are applied.
 \end{remark}
\section{Numerical experiments}\label{section6}~\\
 In this section, we present several numerical experiments in two and three dimensions  to validate the convergence, MBP and energy stability of the LRI methods  developed for the $Q$-tensor model \eqref{1.6a}. To rigorously evaluate the robustness of our proposed LRI methods compared against classical ETD schemes , we conduct a series of numerical experiments using initial data with limited smoothness.  The phase transition process of $+1$ topological defects under Dirichlet boundary conditions is also characterized by analyzing the principal eigenvalues and eigenvectors of the Q-tensor.

    We discretize the spatial domain $\Omega$ with a uniform mesh of size $N$ in each dimension, resulting in a total of $N^d$ grid points. The time step size is denoted by $\tau$, and we compute the numerical solution up to a fixed terminal time $T$.
    Let $D_{h}$ be the central finite difference approximation of the one-dimensional Laplacian operator with periodic boundary conditions as in \cite{ju2015fast}.
 The Laplacian operator can be approximated with the matrix  $D \in \mathbb{R}^{d N \times d N}$ as follows:
 \begin{align}
	 D=\left\{\begin{array}{ll}
	 I \otimes D_{h}+D_{h} \otimes I, & d=2, \\
	 I \otimes I \otimes D_{h}+I \otimes D_{h} \otimes I+D_{h} \otimes I \otimes I, & d=3.
 \end{array}\right. \no
 \end{align}
  Note that  D  is a circulant matrix, thus, we can calculate the product of the matrix exponential and a vector via the fast Fourier transform (FFT).

 For the convergence order, we calculate
 \begin{eqnarray}
	 \rho=\frac{\Norm{Q^{ \tau}-Q^{ \frac{\tau}{2}}}{\mathcal{X}}}{\Norm{Q^{\frac{\tau}{2}}-Q^{ \frac{\tau}{4}}}{\mathcal{X}}}, \nonumber
 \end{eqnarray}
  where $Q^{ \tau}$ denotes the numerical solution acquired at time $T$ with time step size $\tau$.
  Specifically, when $\log_{2}{\rho}\approx 1,2$, the convergence orders are  1 and 2 (cf. \cite{20211223}).
\subsection{Convergence, MBP, and energy dissipation tests.}
 \subsubsection{ Two-dimensional  case}\label{dim2}~\\
We consider the  $Q$-tensor model \eqref{1.6a} in the 2D domain $\Omega=[0,2\pi]^2$ and set the initial condition to be
\begin{align*}
  Q_0(x,y)=\frac{1}{3}(\mathbf{n_0}\mathbf{n_0}^T-\frac{\mathbf{I}}{2}), \quad \mathbf{n_0}=(\text{cos}(\theta),\text{sin}(\theta))^T,
  \end{align*}
   where $\theta=x+y$.
The other parameters are set as $\alpha=-1$, $\beta=0$, $\gamma=2$, and $c=1$. We vary the time step size as $\tau=2^{-k}\tau_{1}$  with $\tau_{1}=2^{-5}$($k=0,1,\ldots,5$). The spatial mesh size is set to  $h=1/128$ and the terminal time is set to $T=1$.

  We first validate the convergence order of the LRI schemes and two ETD schemes.
 In Table \ref{v2_lri_error}-\ref{v2_lri2_error}, the convergence rates produce for both the uniform spectral norm and  $\mathcal{X}$-norm are approximately 1 for LRI1a and LRI1b, and approximately 2 for LRI2a and LRI2b, which are consistent with Theorems \ref{theorem_lri_1a}-\ref{theorem_lri_2}.

 \begin{table}[htbp]
	  \centering
      \footnotesize
		 \begin{tabular}{l *{8}{c}}
	 \toprule
	 \multicolumn{1}{c}{}& \multicolumn{2}{c}{$\mathcal{X}$-norm} & \multicolumn{2}{c}{  Spectral norm}& \multicolumn{2}{c}{$\mathcal{X}$-norm} & \multicolumn{2}{c}{  Spectral norm} \\
	 \cmidrule(lr){2-9}
	 \multicolumn{1}{c}{$\tau_1=2^{-5}$} & Error & Rate & Error & Rate & Error & Rate & Error & Rate \\
	 \cmidrule(lr){2-9}
	 \multicolumn{1}{c}{} &\multicolumn{1}{c}{LRI1a}&\multicolumn{2}{c}{}&\multicolumn{1}{c}{} &\multicolumn{1}{c}{LRI1b}&\multicolumn{3}{c}{} \\
	 \midrule
		   $\tau_1$  &  4.3775E-06&-&3.0953E-06&-&2.6236E-07&-&1.8551E-07&-\\
		   $\tau_1$ / 2  &  2.1792E-06&1.006 &1.5409E-06&1.006 &1.6638E-07&0.657 &1.1765E-07&0.657
			 \\
		   $\tau_1$ / 4  &  1.0873E-06&1.003 &7.6885E-07&1.003 &9.1873E-08&0.857 &6.4964E-08&0.857
			 \\
		   $\tau_1$ / 8  &  5.4310E-07&1.001 &3.8403E-07&1.001 &4.8092E-08&0.934 &3.4006E-08&0.934
			\\
		   $\tau_1$ / 16  &  2.7141E-07&1.001 &1.9191E-07&1.001 &2.4583E-08&0.968 &1.7383E-08&0.968
			 \\
		   $\tau_1$ / 32  &  1.3567E-07&1.000 &9.5933E-08&1.000 &1.2425E-08&0.984 &8.7861E-09&0.984
			 \\
		   $\tau_1$ / 64  &  6.7826E-08&1.000 &4.7960E-08&1.000 &6.2462E-09&0.992 &4.4167E-09&0.992
			 \\
		   $\tau_1$ / 128  & 3.3911E-08&1.000 &2.3979E-08&1.000 &3.1314E-09&0.996 &2.2143E-09&0.996
			 \\
			 \bottomrule
	 \end{tabular}
	 \caption{Errors and convergence rates of the LRI1a and LRI1b schemes.}\label{v2_lri_error}
 \end{table}
 \begin{table}[htbp]
	 \centering
     \footnotesize
	 \begin{tabular}{l *{8}{c}}
		 \toprule
		 \multicolumn{1}{c}{}& \multicolumn{2}{c}{$\mathcal{X}$-norm} & \multicolumn{2}{c}{  Spectral norm}& \multicolumn{2}{c}{$\mathcal{X}$-norm} & \multicolumn{2}{c}{  Spectral norm} \\
		 \cmidrule(lr){2-9}
		 \multicolumn{1}{c}{$\tau_1=2^{-5}$} & Error & Rate & Error & Rate & Error & Rate & Error & Rate \\
		 \cmidrule(lr){2-9}
		 \multicolumn{1}{c}{} &\multicolumn{1}{c}{LRI2a}&\multicolumn{2}{c}{}&\multicolumn{1}{c}{} &\multicolumn{1}{c}{LRI2b}&\multicolumn{3}{c}{} \\
		 \midrule
		 $\tau_1$  &  2.0412E-07&-&1.4434E-07&-&1.2981E-07&-&9.1790E-08&-
		 \\
		 $\tau_1$ / 2  &  5.1895E-08&1.976 &3.6696E-08&1.976 &3.1725E-08&2.033 &2.2433E-08&2.033
		 \\
		 $\tau_1$ / 4  &  1.3079E-08&1.988 &9.2486E-09&1.988 &7.8371E-09&2.017 &5.5417E-09&2.017
		 \\
		 $\tau_1$ / 8  & 3.2829E-09&1.994 &2.3213E-09&1.994 &1.9473E-09&2.009 &1.3770E-09&2.009
		 \\
		 $\tau_1$ / 16  &  8.2233E-10&1.997 &5.8148E-10&1.997 &4.8533E-10&2.004 &3.4318E-10&2.004
		 \\
		 $\tau_1$ / 32  &  2.0578E-10&1.999 &1.4551E-10&1.999 &1.2114E-10&2.002 &8.5661E-11&2.002
		 \\
		 $\tau_1$ / 64  &  5.1471E-11&1.999 &3.6396E-11&1.999 &3.0262E-11&2.001 &2.1399E-11&2.001
		 \\
		 $\tau_1$ / 128  & 1.2871E-11&2.000 &9.1015E-12&2.000 &7.5632E-12&2.000 &5.3480E-12&2.000
		 \\
		 \bottomrule
	 \end{tabular}
	 \caption{Errors and convergence rates of the LRI2a and LRI2b schemes.}\label{v2_lri2_error}
  \end{table}
 \begin{table}[htbp]
	  \centering
		 \begin{tabular}{l *{6}{c}}
	 \toprule
	 \multicolumn{1}{c}{$\tau_1=2^{-5}$} & LRI1a & LRI1b  & ETD1 &    LRI2a & LRI2b& ETD2   \\
	 \midrule
		   $\tau_1$  & 20.70 &20.92 &23.35 &31.26 &32.03& 41.05   \\
		   $\tau_1$ / 2  &  45.18&45.22 &50.35&67.36 &70.87&86.66
			 \\
		   $\tau_1$ / 4  &  90.31&90.38 &100.75 &134.63 &144.71&175.59
			 \\
		   $\tau_1$ / 8  &  180.94&181.01 &201.70&269.26 &291.82&352.87
			\\
		   $\tau_1$ / 16  &  361.95&362.03 &400.04&538.70 &578.54&700.15
			 \\
			 \bottomrule
	 \end{tabular}
	 \caption{Running times (in seconds) of the first and second-order LRI and ETD schemes with N = 2048.} \label{v2_time}
 \end{table}
In Table \ref{v2_time},  we report  the running times (in seconds) of LRI and ETD methods with fixed N = 2048 and varying time step sizes.  We observe that both the first-order LRI1 schemes and second-order LRI2 schemes are faster than  ETD schemes. Especially, the running time of the second-order ETD2 scheme is approximately 1.4 times that of  two second-order LRI schemes.
 Thus, we conclude that  the first and second order LRI schemes are more effective than ETD1 and ETD2 schemes, respectively, in terms of both accuracy and computational cost.

In the following, we test the MBP and energy dissipation properties of the LRI schemes.
We simulate the Q tensor problem with the time step size $\tau=2^{-5}$ over the time interval $[0,100]$.

 \begin{figure}[htbp]
	 \centering
		 \includegraphics[width=0.9\textwidth]{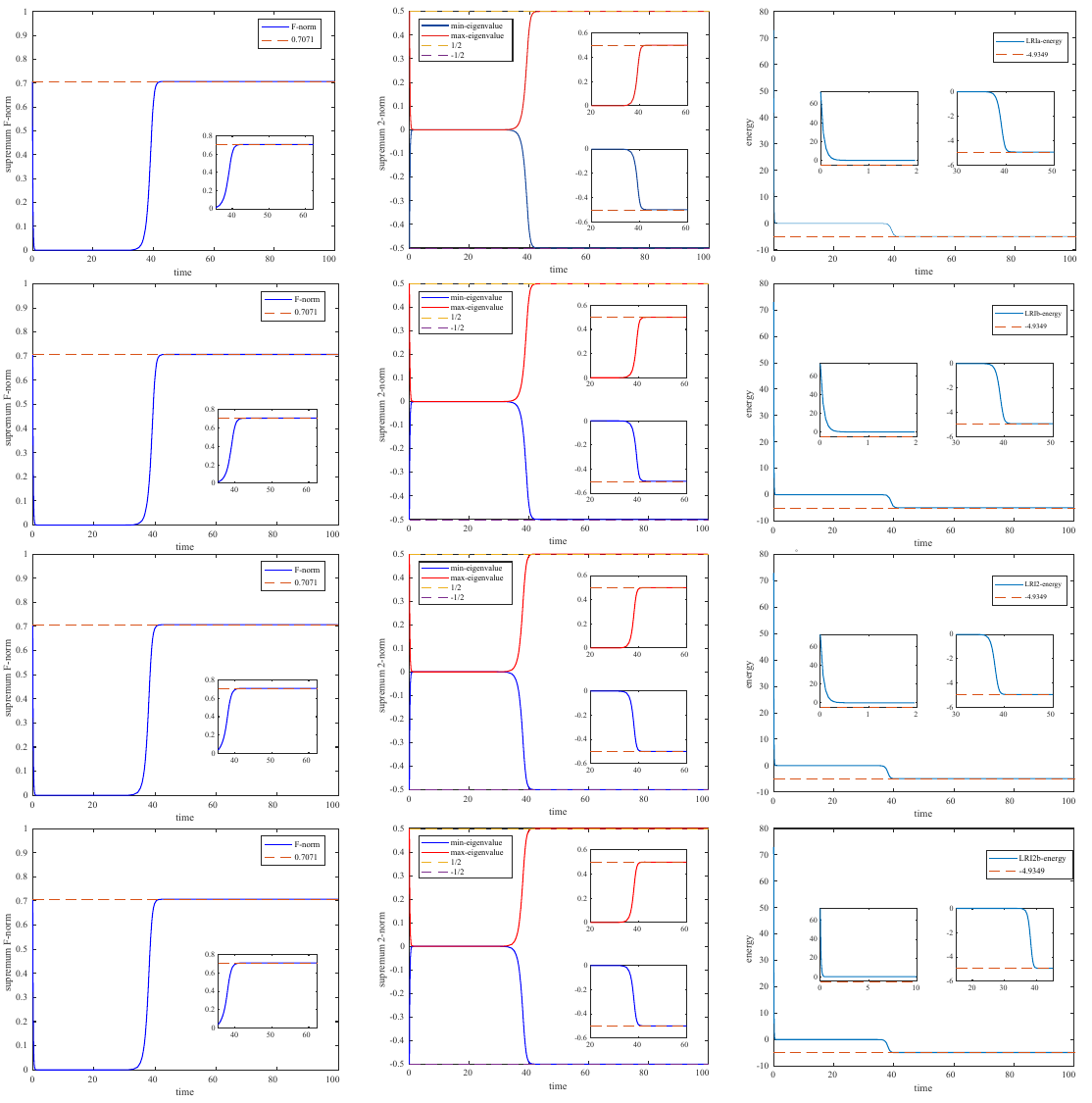}
	 \caption{Two dimension time evolution of the  $\mathcal{X}$-norm, $s_2$-norm, and energy for four LRI schemes  from left to right. The rows, from top to bottom, correspond to the LRI1a, LRI1b, LRI2a, and LRI2b schemes.}
	 \label{v2_lri_1a_max}
 \end{figure}
Figure \ref{v2_lri_1a_max} presents the time evolution of the  $\mathcal{Z}$-norm, uniform spectral norm, and energy for the four LRI schemes.
It can be observed that both the $\mathcal{Z}$-norm and  uniform spectral norm initially decrease but then increase and converge to a steady state over time, which are consistent with the theoretical analysis in Theorem \ref{mbp} and Remark \ref{re3.1}.
Furthermore, the positive and negative eigenvalues of Q  lie strictly within the permissible range $(-\frac{1}{2},\frac{1}{2})$ and stabilize at $ \pm \frac{1}{2}$, which is dictated by the Q-tensor framework. The $\mathcal{Z}$-norm is consistently larger than the uniform spectral norm, which is expected due to the nature of the norms.  The energy  gradually decreases from its initial positive value and eventually stabilizes at a steady negative state which is consistent with Theorem \ref{energy_stability}.
\subsubsection{ Three-dimensional  case}\label{dim3}~\\
We now consider the 3D $Q$-tensor model \eqref{1.6a} on the domain $\Omega=[0,2\pi]^3$  and  set the initial condition to be
 \begin{align*}
 Q_0(x,y,z)=\frac{1}{3}(\frac{\mathbf{n_0}\mathbf{n_0}^T}{\|\mathbf{n_0}\|^2}-\frac{\mathbf{I}}{3}),\quad  \mathbf{n_0}=(\text{cos}(\theta),\text{sin}(\theta),1)^T,
 \end{align*}
where $\theta=x+y+z$.
 The other   parameters are set as $\alpha=-1$, $\beta=1$, $\gamma=2$,   $c=1$. We vary the time step size as $\tau=2^{-k}\tau_{1}$  with $\tau_{1}=2^{-5}$($k=0,1,\ldots,5$). The spatial mesh size is set to  $h=1/64$ and the terminal time is set to $T=1$.

	  In Table \ref{v3_lri1}-\ref{v3_lri2}, the convergence rates produce for both the $\mathcal{X}$-norm and uniform spectral norm are approximately 1 for LRI1a and LRI1b, and approximately 2 for LRI2a and LRI2b, which are consistent with Theorems \ref{theorem_lri_1a}-\ref{theorem_lri_2}.

	 \begin{table}[htbp]
		 \centering
         \footnotesize
		 \begin{tabular}{l *{8}{c}}
			 \toprule
			 \multicolumn{1}{c}{}& \multicolumn{2}{c}{$\mathcal{X}$-norm} & \multicolumn{2}{c}{  Spectral norm}& \multicolumn{2}{c}{$\mathcal{X}$-norm} & \multicolumn{2}{c}{  Spectral norm} \\
			 \cmidrule(lr){2-9}
			 \multicolumn{1}{c}{$\tau_1=2^{-5}$} & Error & Rate & Error & Rate & Error & Rate & Error & Rate \\
			 \cmidrule(lr){2-9}
			 \multicolumn{1}{c}{} &\multicolumn{1}{c}{LRI1a}&\multicolumn{2}{c}{}&\multicolumn{1}{c}{} &\multicolumn{1}{c}{LRI1b}&\multicolumn{3}{c}{} \\
			 \midrule
			 $\tau_1$ & 2.5137E-03&-&2.0524E-03&-&1.8065E-03&-&1.4749E-03&-
			  \\
			 $\tau_1$ / 2 & 1.1866E-03&1.083 &9.6883E-04&1.083 &9.2297E-04&0.969 &7.5358E-04&0.969
			 \\
			 $\tau_1$ / 4 & 5.7637E-04&1.042 &4.7060E-04&1.042 &4.6606E-04&0.986 &3.8053E-04&0.986
			 \\
			 $\tau_1$ / 8 & 2.8404E-04&1.021 &2.3191E-04&1.021 &2.3413E-04&0.993 &1.9116E-04&0.993
			  \\
			 $\tau_1$ / 16 & 1.4100E-04&1.010 &1.1512E-04&1.010 &1.1733E-04&0.997 &9.5799E-05&0.997
			 \\
			 $\tau_1$ / 32 & 7.0243E-05&1.005 &5.7352E-05&1.005 &5.8732E-05&0.998 &4.7954E-05&0.998
			 \\
			 $\tau_1$ / 64 & 3.5058E-05&1.003 &2.8624E-05&1.003 &2.9383E-05&0.999 &2.3990E-05&0.999
			  \\
			 $\tau_1$ / 128 & 1.7513E-05&1.001 &1.4299E-05&1.001 &1.4695E-05&1.000 &1.1998E-05&1.000
			 \\
			 \bottomrule
		 \end{tabular}
			 \caption{Errors and convergence rates of the LRI1a and LRI1b schemes.}
			 \label{v3_lri1}
	 \end{table}
	 \begin{table}[htbp]
	 \centering
     \footnotesize
	 \begin{tabular}{l *{8}{c}}
		 \toprule
		 \multicolumn{1}{c}{}& \multicolumn{2}{c}{$\mathcal{X}$-norm} & \multicolumn{2}{c}{  Spectral norm}& \multicolumn{2}{c}{$\mathcal{X}$-norm} & \multicolumn{2}{c}{  Spectral norm} \\
		 \cmidrule(lr){2-9}
		 \multicolumn{1}{c}{$\tau_1=2^{-5}$} & Error & Rate & Error & Rate & Error & Rate & Error & Rate \\
		 \cmidrule(lr){2-9}
		 \multicolumn{1}{c}{} &\multicolumn{1}{c}{LRI2a}&\multicolumn{2}{c}{}&\multicolumn{1}{c}{} &\multicolumn{1}{c}{LRI2b}&\multicolumn{3}{c}{} \\
		 \midrule
		 $\tau_1$ &3.3744E-04&-&2.7552E-04&-& 1.5210E-04&-&1.2419E-04&-
		 \\
		 $\tau_1$ / 2  &8.1494E-05&2.050 &6.6539E-05&2.050 & 4.1446E-05&1.876 &3.3841E-05&1.876
		 \\
		 $\tau_1$ / 4  &1.9840E-05&2.038 &1.6199E-05&2.038 & 1.0765E-05&1.945 &8.7896E-06&1.945
		 \\
		 $\tau_1$ / 8  &4.8831E-06&2.023 &3.9871E-06&2.023 & 2.7398E-06&1.974 &2.2371E-06&1.974
		 \\
		 $\tau_1$ / 16  &1.2106E-06&2.012 &9.8843E-07&2.012 & 6.9089E-07&1.988 &5.6412E-07&1.988
		 \\
		 $\tau_1$ / 32  &3.0133E-07&2.006 &2.4603E-07&2.006 & 1.7346E-07&1.994 &1.4163E-07&1.994
		 \\
		 $\tau_1$ / 64  &7.5166E-08&2.003 &6.1373E-08&2.003 & 4.3456E-08&1.997 &3.5483E-08&1.997
		 \\
		 $\tau_1$ / 128  &1.8771E-08&2.002 &1.5326E-08&2.002 & 1.0875E-08&1.998 &8.8799E-09&1.998
		 \\
		 \bottomrule
	 \end{tabular}
	 \caption{Errors and convergence rates of the LRI2a and LRI2b schemes.}
	 \label{v3_lri2}
 \end{table}
 \begin{table}[htbp]
	  \centering
		 \begin{tabular}{l *{6}{c}}
	 \toprule
	 \multicolumn{1}{c}{$\tau_1=2^{-5}$} & LRI1a & LRI1b  & ETD1 &    LRI2a & LRI2b& ETD2   \\
	 \midrule
		   $\tau_1$  & 55.54 & 55.56 &72.69 &84.21 &92.01& 120.30  \\
		   $\tau_1$ / 2  &  107.03&107.42 &144.32&173.75 &192.07&246.60
			 \\
		   $\tau_1$ / 4  &  214.66&215.15 &290.54 &339.63 &370.85&478.12
			 \\
		   $\tau_1$ / 8  &  423.96&425.08 &572.67&646.13 &708.31&898.23
			\\
		   $\tau_1$ / 16  &  849.80&851.04 &1074.13&1212.80 &1335.59&1645.69
			 \\
			 \bottomrule
	 \end{tabular}
	 \caption{Running times (in seconds) of the first and second-order LRI and ETD schemes with h = 1/128.} \label{v3_time}
 \end{table}
In Table \ref{v3_time},  we report  the running times  of LRI and ETD schemes with fixed h = 1/128 and varying time step sizes.  It can be  observed that both the first-order  and second-order LRI schemes are faster than the  ETD schemes.

In the following, we test the MBP and energy dissipation properties of the LRI schemes  with the time step size $\tau=2^{-5}$ over the time interval $[0,100]$ and the same initial conditions and parameters.

 \begin{figure}[htbp]
		 \includegraphics[width=0.98\textwidth]{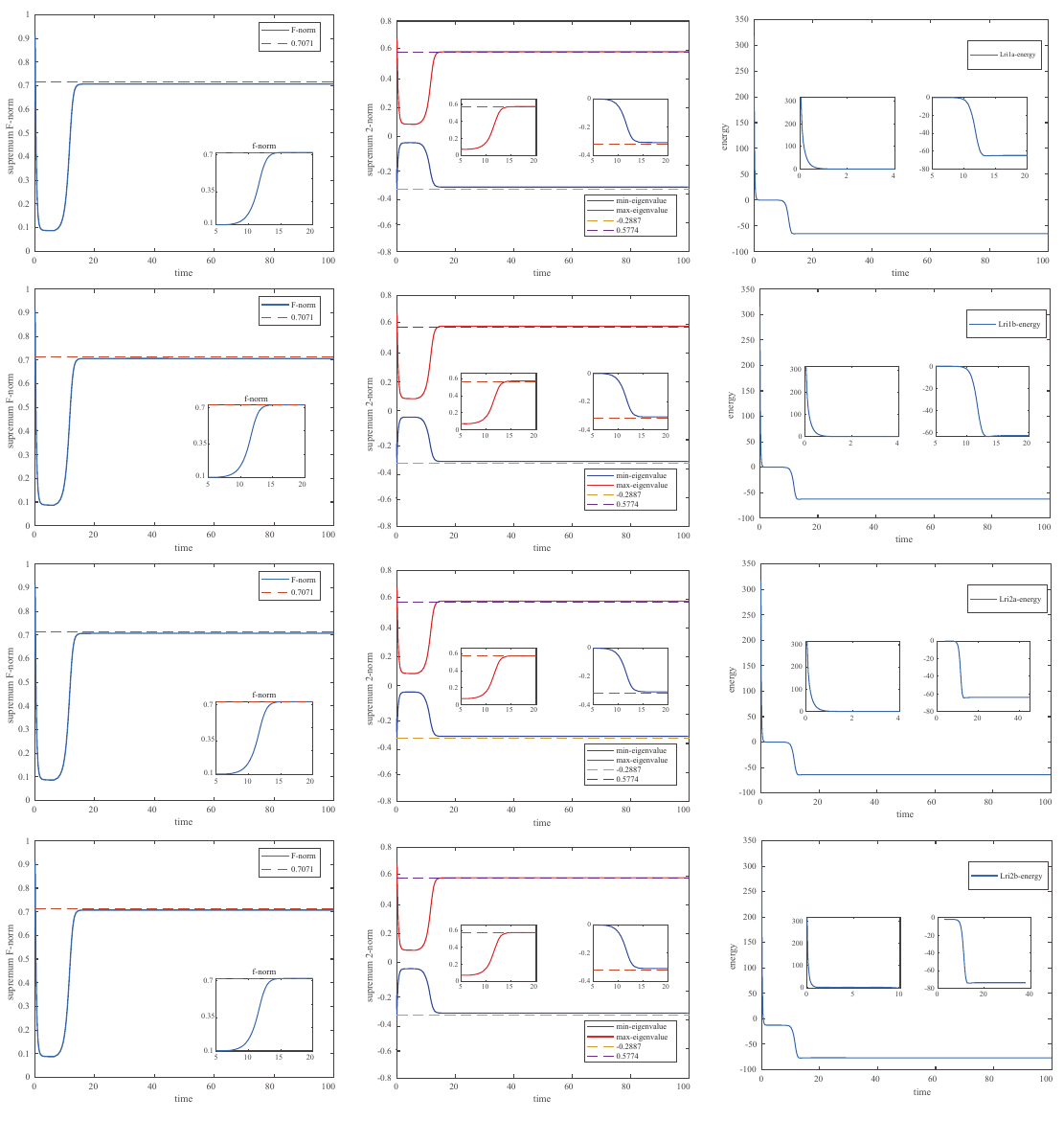}
	 \caption{Three dimension time evolution of the  $\mathcal{X}$-norm, $s_2$-norm, and energy for four LRI schemes  from left to right. The rows, from top to bottom, correspond to the LRI1a, LRI1b, LRI2a, and LRI2b schemes.}
	 \label{v3_lri1a}
 \end{figure}

Figure \ref{v3_lri1a} presents the evolution trends of the $\mathcal{Z}$-norm, the uniform spectral norm, and energies of the solutions of four LRI schemes.
The eigenvalue visualization results distinctly illustrate that the maximum eigenvalue $0.5774$ is twice the absolute value of the minimum eigenvalue $-0.2887$, aligning well with that  the $Q$-tensor should  have two degenerate eigenvalues in \cite{mottramintroduction}.
 It can be   observed that both the $\mathcal{Z}$-norm and uniform spectral norm   stabilize at equilibrium and the $\mathcal{Z}$-norm is consistently larger than the uniform spectral norm, which is consistent with Theorem \ref{mbp}. Furthermore, driven by the physical constraints imposed on the Q-tensor parameter, the positive and negative eigenvalues lie strictly within the permissible range $(-\frac{1}{3},\frac{2}{3})$ dictated by the Q-tensor framework.
   The energy  decreases from its initial positive value and eventually stabilizes at a steady negative state, which is consistent with Theorem \ref{energy_stability}.
    \subsection{Convergence Tests with Low-Regularity Initial Data}~\\
 In many physical scenarios, initial configurations  contain defects or sharp phase boundaries, which possess limited regularity. Consequently, the ability of a numerical scheme to preserve its theoretical accuracy under such conditions serves as a critical benchmark of its reliability. To this end, we demonstrate  the robustness of the LRI schemes  by comparing their performance with classical ETD schemes.

We employ a standard technique to construct initial data with precisely controlled Sobolev regularity.
Specifically, we define the scalar function as follows:
\begin{equation}
    S(\mathbf{x}) = \left( \max(0, R - \|\mathbf{x}\|) \right)^k, \quad \mathbf{x} \in \Omega \subset \mathbb{R}^3,
    \label{eq:scalar_initial}
\end{equation}
where $R > 0$ is a constant defining the radius of the support, and  $S(\mathbf{x})$  belongs to the Sobolev space $H^k(\Omega)$ but not to $H^{k+1}(\Omega)$.
Based on this scalar field, the initial Q-tensor field $Q_0(\mathbf{x})$ is constructed as:
\begin{equation}
    Q_0(\mathbf{x}) = S(\mathbf{x}) \left( \mathbf{n} \otimes \mathbf{n} - \frac{1}{3}\mathbf{I} \right),
    \label{eq:tensor_initial}
\end{equation}
where $\mathbf{n}$ is a fixed unit vector, and $\mathbf{I}$ is the $3 \times 3$ identity tensor. This construction ensures that $Q_0(\mathbf{x})$ is symmetric and traceless, and its Sobolev regularity is determined entirely by  the scalar field $S(\mathbf{x})$.

    The other parameters are set as  $\alpha=-1$, $\beta=1$, $\gamma=2$. We vary the time step size as $\tau=2^{-k}\tau_{1}$  with $\tau_{1}=2^{-5}$($k=0,1,\ldots,5$). The spatial mesh size is set to be $N=64$.
    To investigate the performance of the numerical schemes, we consider two representative initial conditions generated from \eqref{eq:scalar_initial} by setting $k=1$ and $k=3$, which belong to $H^1(\Omega) \setminus H^2(\Omega)$ and $H^3(\Omega) \setminus H^4(\Omega)$, respectively.
\begin{figure}[htbp]
    \centering 

    \begin{minipage}[b]{0.48\textwidth}
        \includegraphics[width=\textwidth]{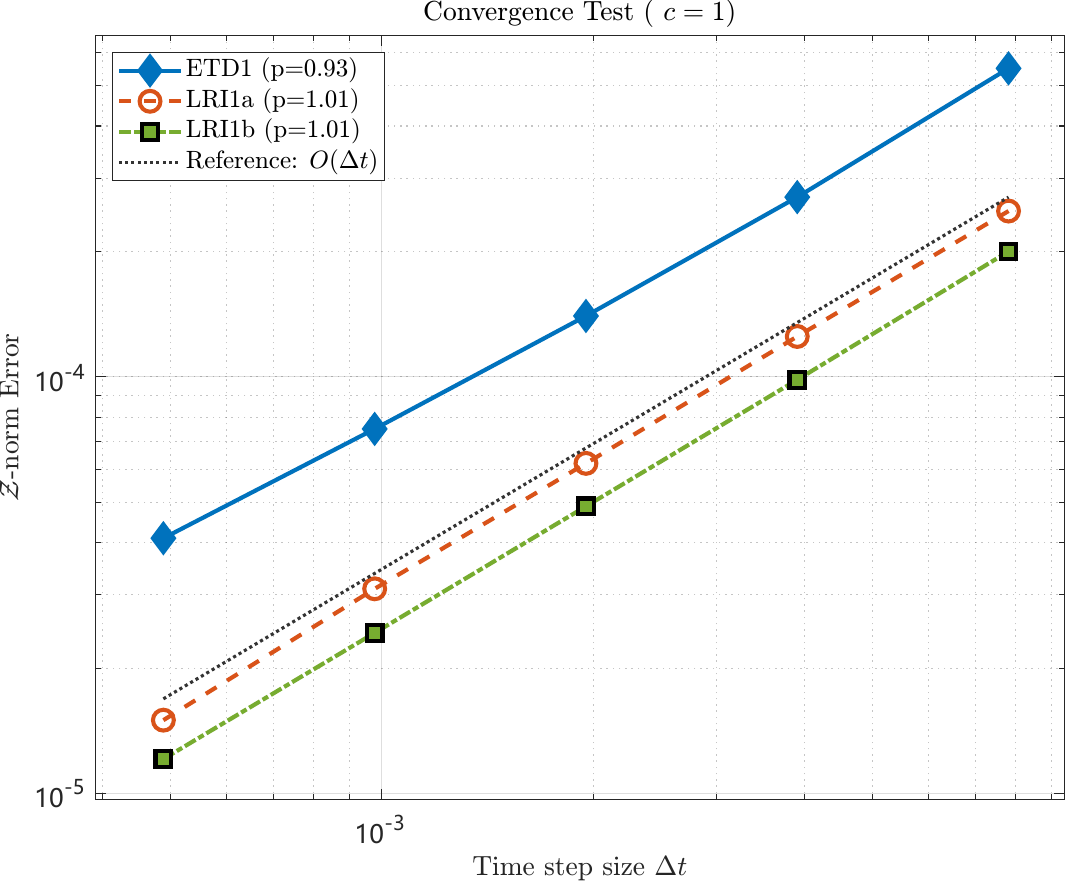}

    \end{minipage}
    \hfill
    \begin{minipage}[b]{0.48\textwidth}
        \includegraphics[width=\textwidth]{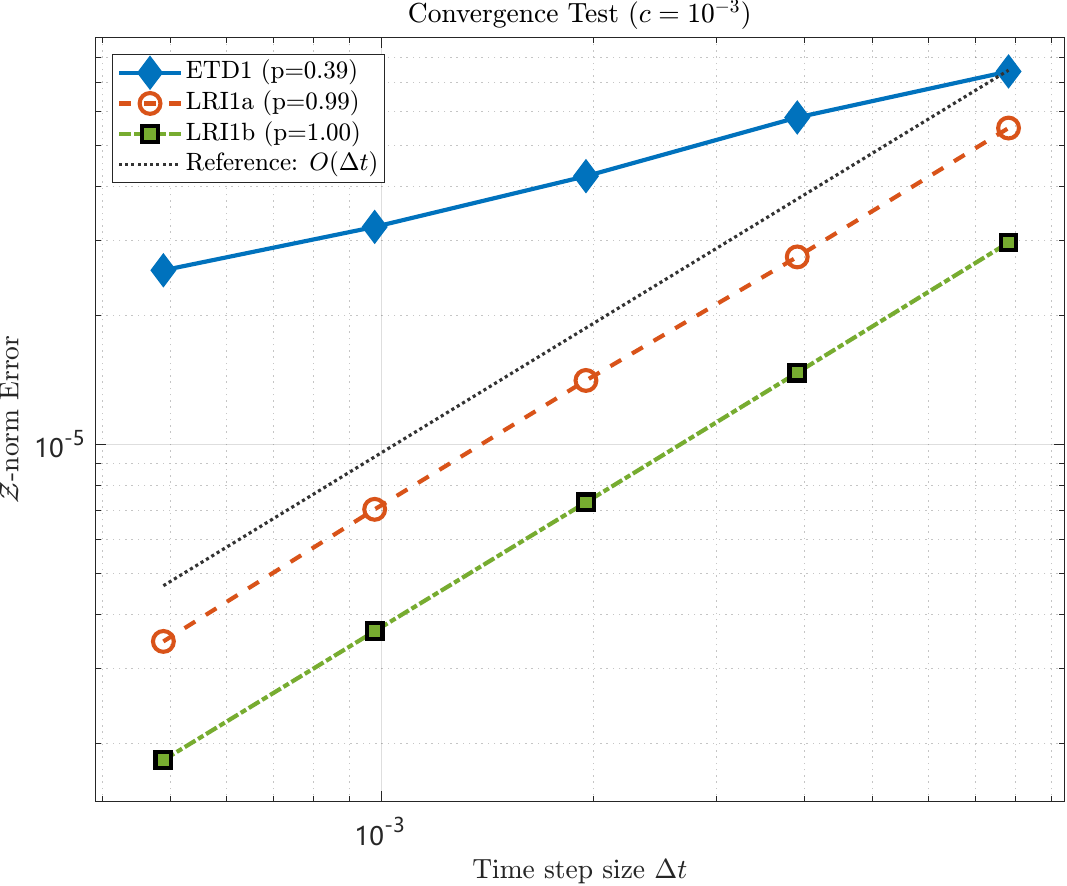}
    \end{minipage}

    \caption{Convergence orders in the $\mathcal{Z}$-norm for the LRI1a, LRI1b, and ETD1 schemes with $H^1 \setminus H^2$ initial data, evaluated at $T=1$ for $c = 1$ and $c = 10^{-3}$.}
    \label{fig_convergence_comparison}
\end{figure}
\begin{figure}[htbp] 
    \centering 

    \begin{minipage}[b]{0.48\textwidth} 
        \includegraphics[width=\textwidth]{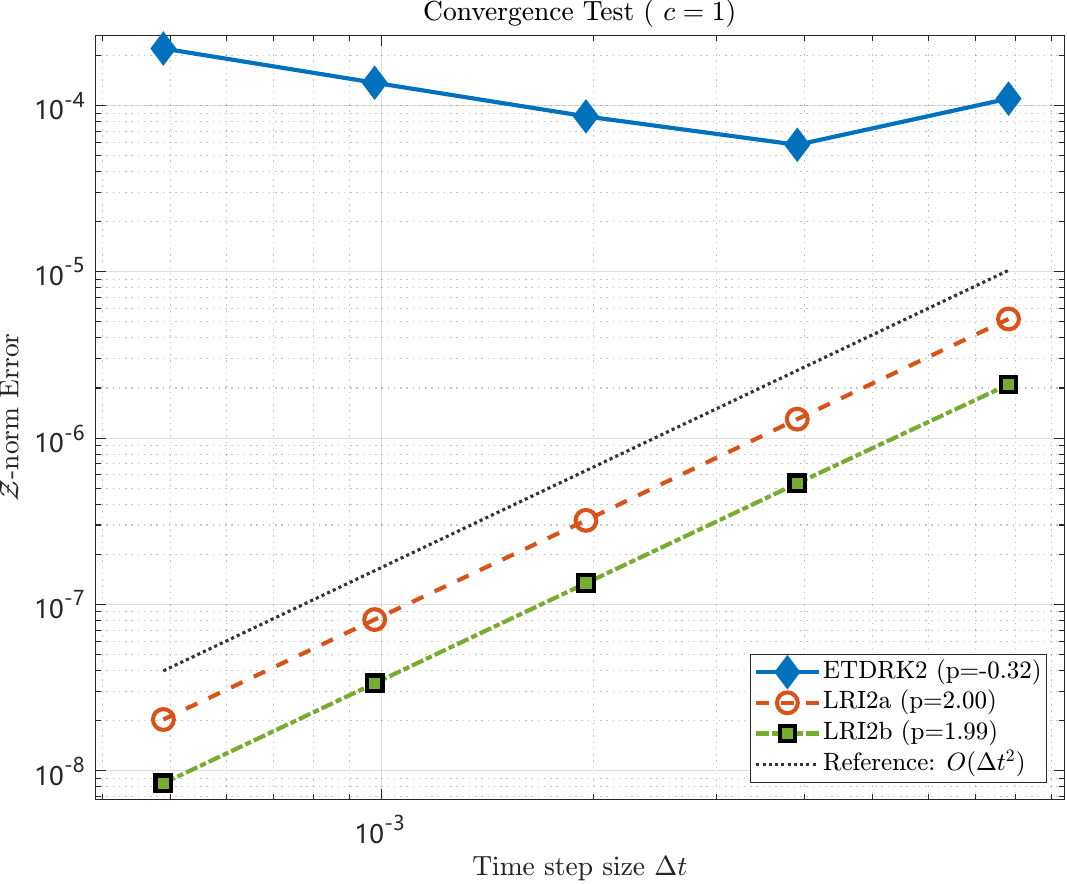}
    \end{minipage}
    \hfill 
    \begin{minipage}[b]{0.48\textwidth} 
        \includegraphics[width=\textwidth]{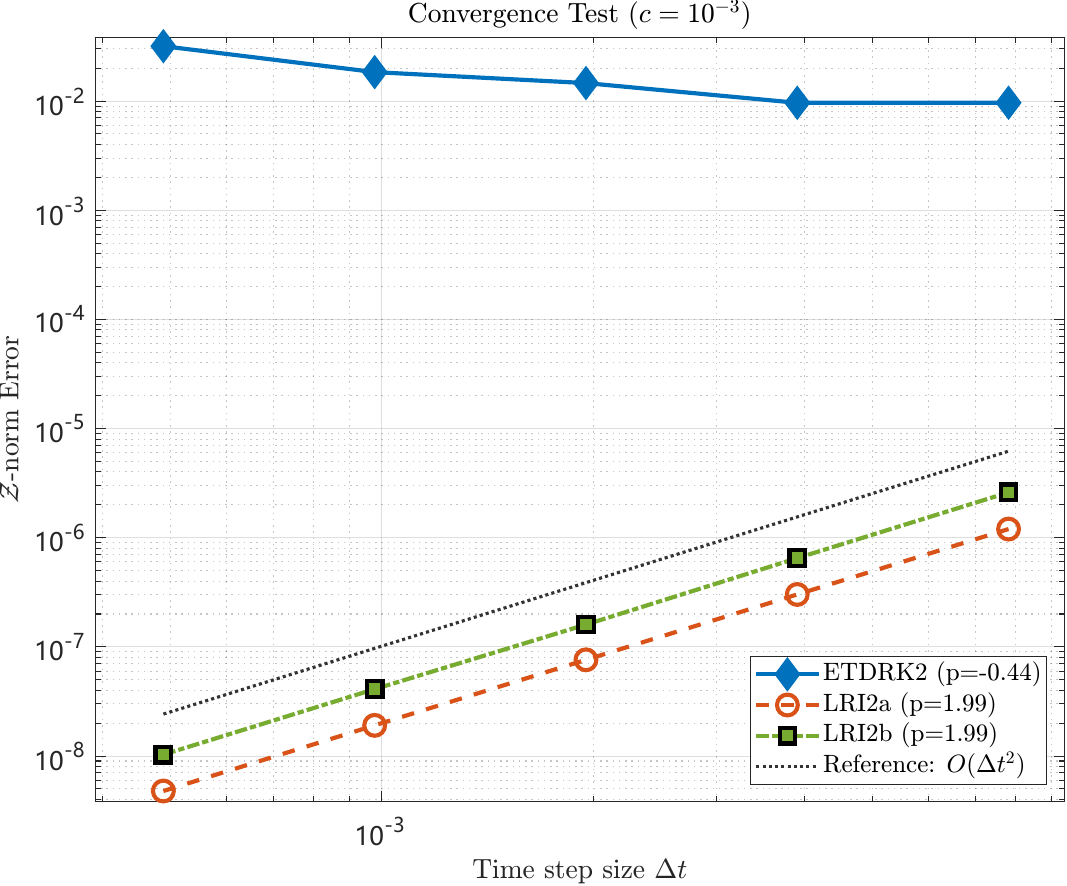}
    \end{minipage}

    \caption{Convergence orders in the $\mathcal{Z}$-norm for the LRI2a, LRI2b, and ETDRK2 schemes with $H^3 \setminus H^4$ initial data, evaluated at $T=1$ for $c = 1$ and $c = 10^{-3}$. }
    \label{fig_convergence_comparison2}
\end{figure}

Figures \ref{fig_convergence_comparison} and \ref{fig_convergence_comparison2}  illustrate the convergence orders in the $\mathcal{Z}$-norm for the LRI and ETD schemes. The results  demonstrate that both LRI schemes consistently achieve their theoretical convergence orders of 1 and 2, respectively, regardless of the regularity of the initial data. In contrast, the ETD schemes exhibit significant order reduction when applied to the lower-regularity initial data. This discrepancy is particularly pronounced in the nonlinearity-dominated regime (i.e., $c = 10^{-3}$), where the ETD1 scheme's convergence order drops below 0.5 for $H^1 \setminus H^2$ initial data, and the ETDRK2 scheme's order falls below 1 for $H^3 \setminus H^4$ initial data.
\subsection{Dynamics of defects in Liquid Crystals}~\\
 In this subsection, we  simulate the 2D and 3D dynamics of defects in nematic LCs via the  LRI2a scheme.   Inspired by \cite{qiaozhonghua,10.4208/eajam.291018.070619}, we investigate   the formation and evolution of the defect patterns from specific initial and boundary conditions.
\subsubsection{Two dimensional defect dynamics}~\\
 We first simulate the  Q-tensor gradient flow model \eqref{1.6a} in a 2D rectangular domain $[0,L_x] \times [0,L_y]$, where $L_x = L_y = 2 \pi$.   The initial and dirichlet boundary condition are set as
\[
Q|_{\partial \Omega} = \frac{\mathbf{n}_0 \mathbf{n}_0^T}{\|\mathbf{n}_0\|^2} - \frac{\mathbf{I}}{2},\quad
Q|_{ \Omega} = \frac{\mathbf{n}_1 \mathbf{n}_1^T}{\|\mathbf{n}_1\|^2} - \frac{\mathbf{I}}{2},
\]
where $\mathbf{n}_0 = (x - 0.5 L_x, y - 0.5 L_y)^T, \mathbf{n}_1 =  (x- 0.25 L_x, y - 0.25 L_y)^T$. The other parameters are set as $\alpha=-0.2$, $\beta=0$, $\gamma=0.5$, and $c=0.1$.

 \begin{figure}[htbp]
	\includegraphics[width=\textwidth]{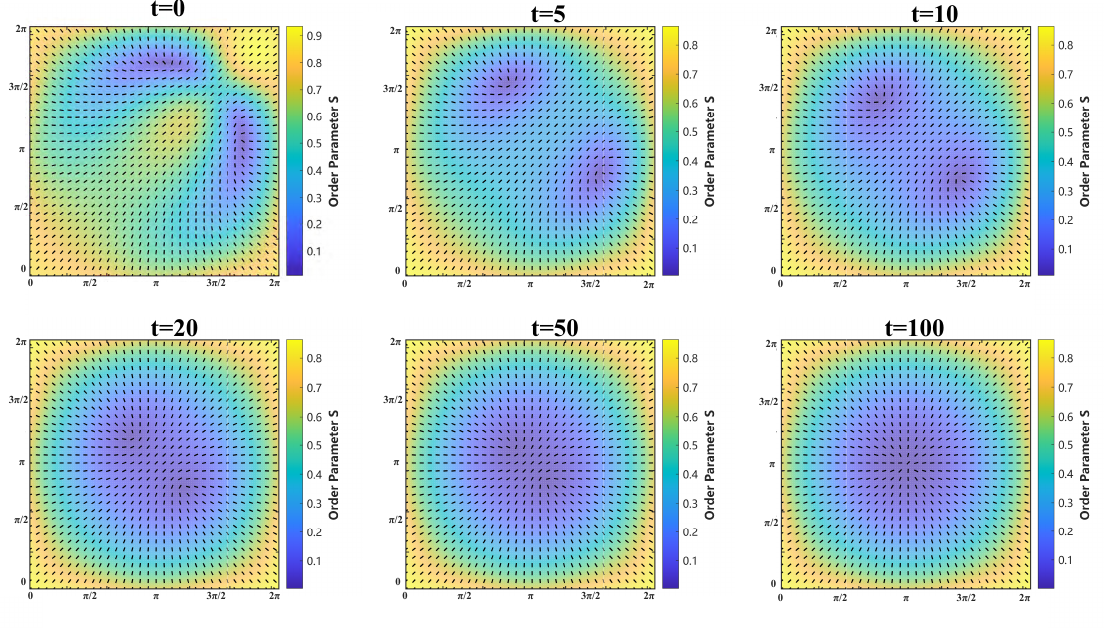}
	 \caption{Evolution of the scalar order parameter $S$ and the director field at $t = 0,5,10,20,50,100$.}
	 \label{v2_lri_principle_eigenvalue}
 \end{figure}
 \begin{figure}
	\centering
	\includegraphics[width=0.45\textwidth]{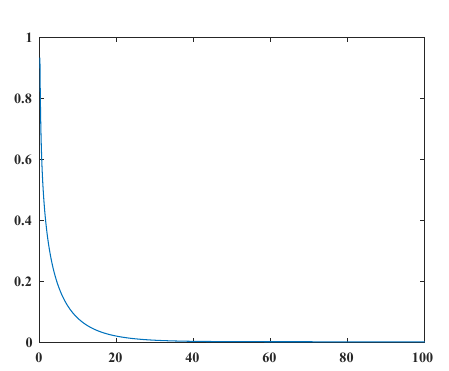}
	 \caption{Evolution trend of the energy of $Q$.}
	 \label{energy}
 \end{figure}
 Figure \ref{v2_lri_principle_eigenvalue} shows the evolution of the scalar order parameter $S$ and the director field. At $t=0$, two $+\frac{1}{2}$ defects are observed. Over time, these defects move towards each other and merge into a single $+1$ defect at the domain center. This coalescence minimizes the total free energy, consistent with the boundary conditions and nematic liquid crystal dynamics. As shown in Figure \ref{energy}, the system's energy decreases monotonically during this process, indicating a relaxation towards equilibrium.
\subsubsection{Three dimensional defect dynamics}~\\
We now simulate the  Q-tensor gradient flow model \eqref{1.6a} in a 3D domain $[0,L_x] \times [0,L_y]\times [0,L_z]$, where $L_x = L_y =L_z= 2\pi$.
  The  initial and dirichlet boundary condition are set as
  \begin{align}
    Q|_{\partial \Omega} = \frac{\mathbf{n}_0 \mathbf{n}_0^T}{\|\mathbf{n}_0\|^2} - \frac{\mathbf{I}}{2},\quad
Q|_{ \Omega} = \frac{\mathbf{n}_1 \mathbf{n}_1^T}{\|\mathbf{n}_1\|^2} - \frac{\mathbf{I}}{2},\label{in1}
  \end{align}
where $\mathbf{n}_0 = (x - 0.5 L_x, y - 0.5 L_y,0)^T, \mathbf{n}_1 =  (x- 0.5 L_x, y - 0.5 L_y,0)^T$. The other parameters are set as $\alpha=-1$, $\beta=1$, $\gamma=2$, and $c=1$.
 To visualize the biaxiality property, we follow \cite{wang2021modelling} and define
 \begin{equation*}
	 \beta_b=1-6 \frac{\left(\operatorname{tr} Q^{3}\right)^{2}}{\left(\operatorname{tr} Q^{2}\right)^{3}},
 \end{equation*}
 where $\beta_b$ is the biaxiality parameter of the Q-tensor  quantifying the degree of biaxiality in the system, with values ranging from 0 (uniaxial) to 1 (biaxial).

 \begin{figure}[htbp]
	\includegraphics[width=\textwidth]{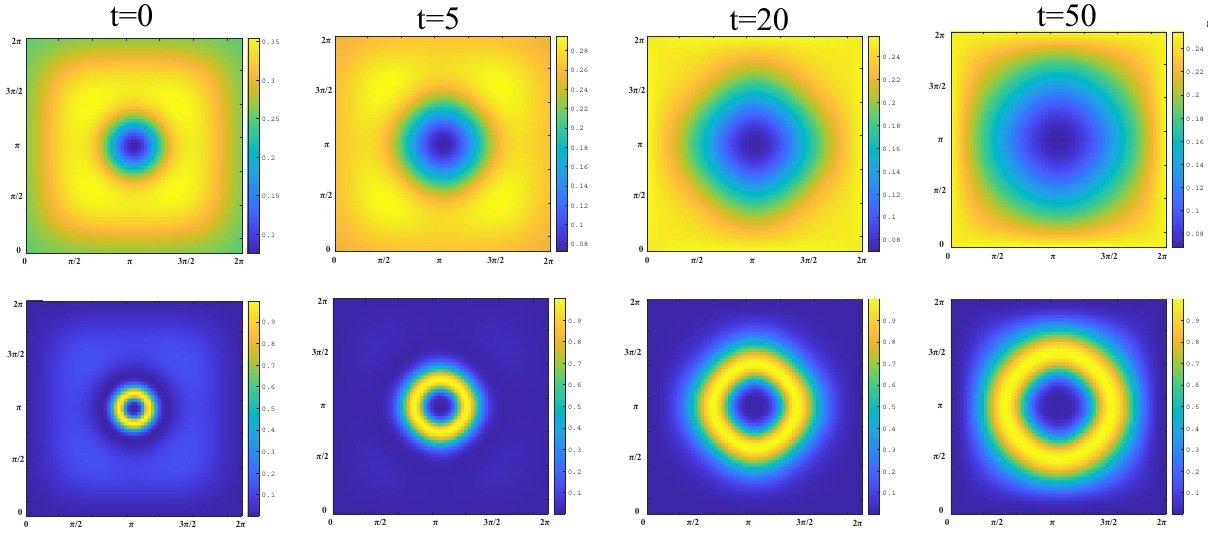}
		 \caption{Evolution trends of the cross-sections of the order parameter and the biaxiality of $Q$ at t = 0,5,20,50.}
	 \label{v3_lri_biaxiality}
 \end{figure}

Figure \ref{v3_lri_biaxiality} illustrates the evolution of the cross-sections of the scalar order parameter (top row) and the biaxiality parameter (bottom row) at t = 0, 5, 20, and 50. As the system evolves, the bulk region develops a higher degree of nematic order, while the disorder becomes confined to a sharp defect core. Concurrently, a distinct ring of high biaxiality forms around this core, indicating  that the system avoids the unphysical singularity at the defect core by developing a finite biaxial state.
   \begin{figure}[htbp]
	\includegraphics[width=\textwidth]{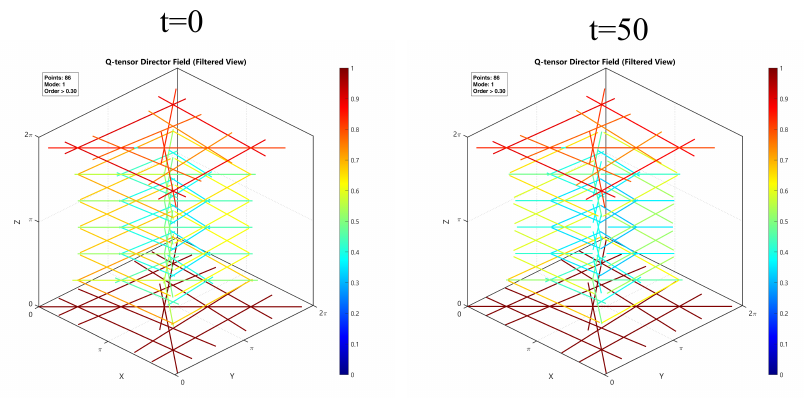}
		 \caption{Evolution trends of the director field  of $Q$ at t = 0,50.}
	 \label{v3_dual}
 \end{figure}

Figure \ref{v3_dual} shows the 3D director field filtered for $S > 0.30$.
The evolution from $t=0$ to $t=50$ reveals the emergence of a central region where the order parameter falls below the threshold ($S < 0.30$).
This visualization highlights the formation of a disordered, isotropic core at the center of the defect structure, a mechanism characteristic of Landau-de Gennes dynamics to regularize energy singularities.

In the following, we simulate the dynamics of defects in LCs with a different initial condition and  boundary condition, which are set as
\begin{align}
    Q|_{\partial \Omega} = \frac{\mathbf{n}_0 \mathbf{n}_0^T}{\|\mathbf{n}_0\|^2} - \frac{\mathbf{I}}{2},\quad
Q|_{ \Omega} = \frac{\mathbf{n}_1 \mathbf{n}_1^T}{\|\mathbf{n}_1\|^2} - \frac{\mathbf{I}}{2},\label{in2}
\end{align}
where $\mathbf{n}_0 = (x - 0.5 L_x, y - 0.5 L_y,0)^T, \mathbf{n}_1 =  (x- 0.4 L_x, y - 0.4 L_y,0)^T$. The other parameters are set as $\alpha=-1$, $\beta=1$, $\gamma=2$, and $c=1$.


    \begin{figure}[htbp]
	\includegraphics[width=\textwidth]{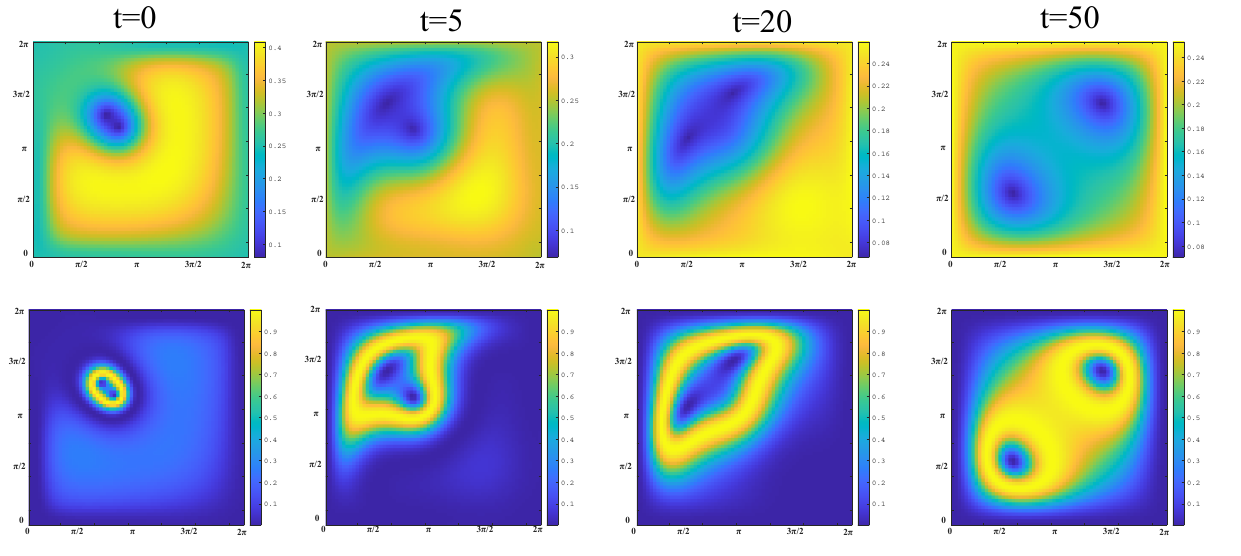}
		 \caption{Evolution trends of the cross-sections of the order parameter and the biaxiality of $Q$ at t = 0,5,20,50.}
	 \label{v3_tensor}
 \end{figure}

    Figure \ref{v3_tensor} show the evolution trends of the cross-sections of the order parameter and the biaxiality of $Q$, illustrating a classic topological relaxation pathway: the fission of an unstable line defect into a pair of point defects. The top row shows the scalar order parameter, where an initial, single low-order region progressively deforms and splits into two separate cores by t=50. The bottom row, showing the biaxiality, reveals the underlying structural changes. The initial defect is sheathed in a single biaxial ring that mirrors the deformation of the core, eventually pinching off and separating into two distinct rings. This process demonstrates the system's drive to minimize its elastic energy by breaking an unstable, higher-energy configuration into a pair of stable, lower-energy topological point defects.

    \begin{figure}[htbp]
	\includegraphics[width=\textwidth]{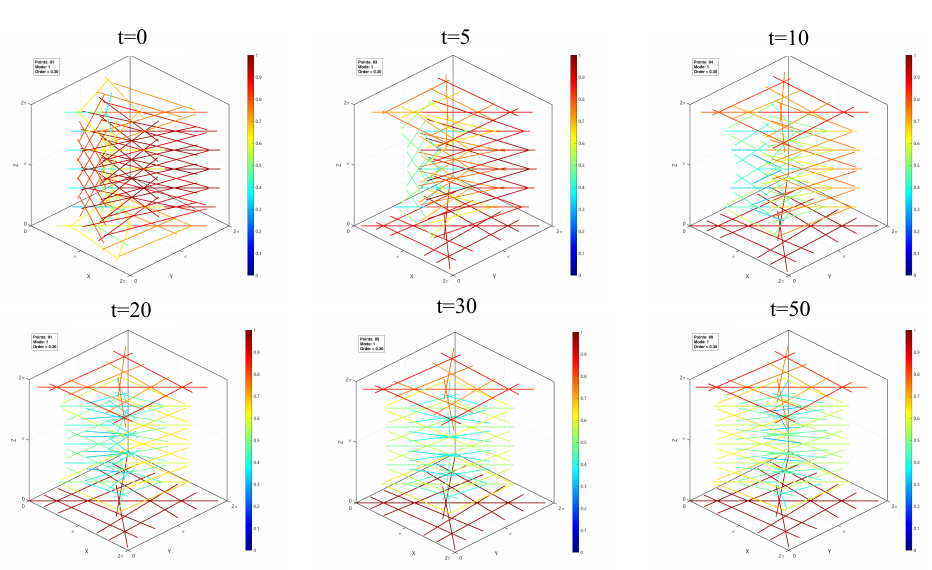}
		 \caption{Evolution trends of the director field  of $Q$ at t = 0,20,50.}
	 \label{v3_deli_qui1}
 \end{figure}
Figure \ref{v3_deli_qui1} depicts the topological relaxation of a defect loop from an energetically unstable state.
At $t=0$, the loop is initialized with a complex, tilted geometry.
As the system minimizes its free energy ($t=5$ to $t=30$), the loop untwists and rotates, eliminating the initial geometric distortion.
By $t=50$, the defect stabilizes into a symmetric loop perpendicular to the $xy$-plane, achieving a local energy minimum.

    \begin{figure}[htbp]
	\includegraphics[width=\textwidth]{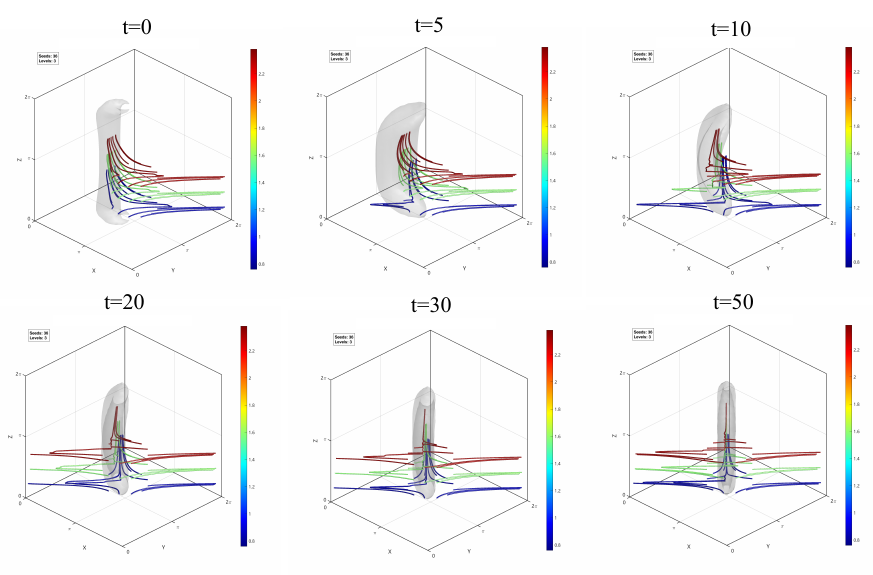}
		 \caption{Evolution trends of the isosurface and streamtubes  of $Q$ at t = 0,5,10,20,50.}
	 \label{v3_deli_qui}
 \end{figure}
Figure \ref{v3_deli_qui} illustrates the dynamic relaxation pathway towards a stable escaped hedgehog defect configuration.
Streamtubes visualize the director field orientation, while an isosurface highlights the region of high biaxiality.
Starting from an asymmetric initial state at $t=0$, the system undergoes a clear symmetrization process between $t=5$ and $t=30$.
During this evolution, the focal point of the director escape migrates to the central axis, and the biaxial core transforms from a contorted structure into a symmetric, spindle-like shape.
The final state at $t=50$ represents the classic escaped hedgehog configuration, demonstrating how the biaxial core effectively regularizes the point singularity to allow for a smooth director field transition.
    \begin{figure}[htbp]
	\includegraphics[width=\textwidth]{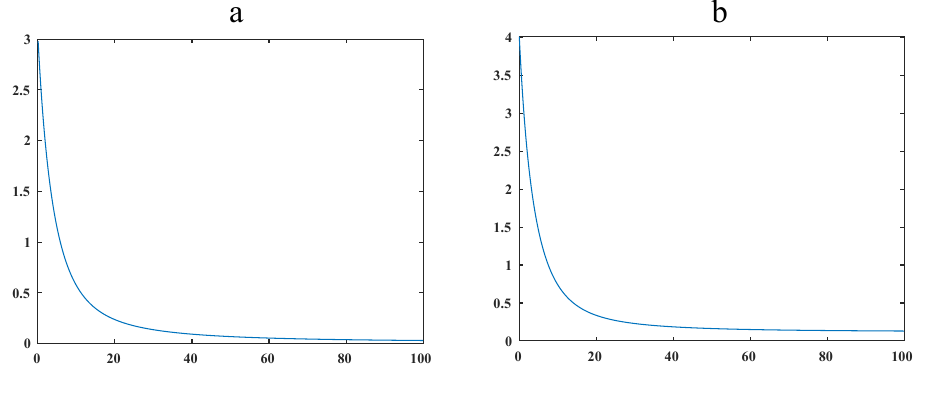}
		 \caption{Evolution of the total free energy over time with the initial condition \eqref{in1} (a) and \eqref{in2} (b).}
	 \label{v3_energy}
 \end{figure}

Figure \ref{v3_energy} demonstrates how the compatibility between initial and boundary conditions affects energy dissipation.
While both cases show monotonic energy decay, system (b) undergoes a more significant structural rearrangement to accommodate the boundary constraints, resulting in a larger energy drop.
Ultimately, (b) stabilizes at a higher energy level than (a), suggesting that the system has settled into a metastable local minimum rather than the global ground state.

In the following, we simulate the dynamics of defects in LCs with a different initial condition and  boundary condition, which are set as
\begin{align}
    Q|_{\partial \Omega} = \frac{\mathbf{n}_0 \mathbf{n}_0^T}{\|\mathbf{n}_0\|^2} - \frac{\mathbf{I}}{2},\quad
Q|_{ \Omega} = \frac{\mathbf{n}_1 \mathbf{n}_1^T}{\|\mathbf{n}_1\|^2} - \frac{\mathbf{I}}{2},\label{in3}
\end{align}
where $\mathbf{n}_0 = (x - 0.5 L_x, y - 0.5 L_y,z-0.5 L_z)^T, \mathbf{n}_1 =  (x- 0.4 L_x, y - 0.4 L_y,0)^T$. The other parameters are set as $\alpha=-1$, $\beta=1$, $\gamma=2$, and $c=1$.

    \begin{figure}[htbp]
	\includegraphics[width=\textwidth]{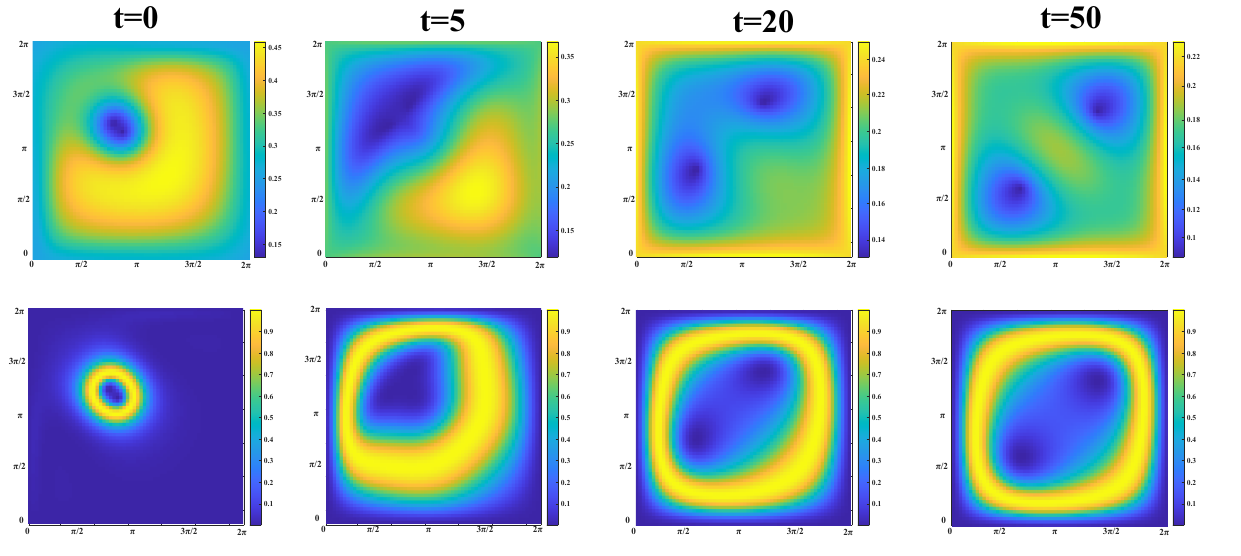}
		 \caption{Evolution trends of the cross-sections of the order parameter and the biaxiality of $Q$ at t = 0,5,10,20,30,50.}
	 \label{v3_tensor1}
 \end{figure}

Figure \ref{v3_tensor1} presents the temporal evolution of the scalar order parameter $S$ (top row) and the biaxiality $\beta^2$ (bottom row) on a 2D cross-section.
Initially ($t=0$), the system exhibits a single singular core, characterized by a minimum in $S$ and surrounded by a ring of high biaxiality.
As the relaxation proceeds ($t=5$ to $t=50$), this central defect elongates and undergoes bifurcation.
The high-biaxiality ring deforms and eventually pinches off, resulting in two distinct defect cores.
This process confirms the topological decomposition of an unstable $+1$ disclination into a pair of stable $+1/2$ defect lines.

    \begin{figure}[htbp]
	\includegraphics[width=\textwidth]{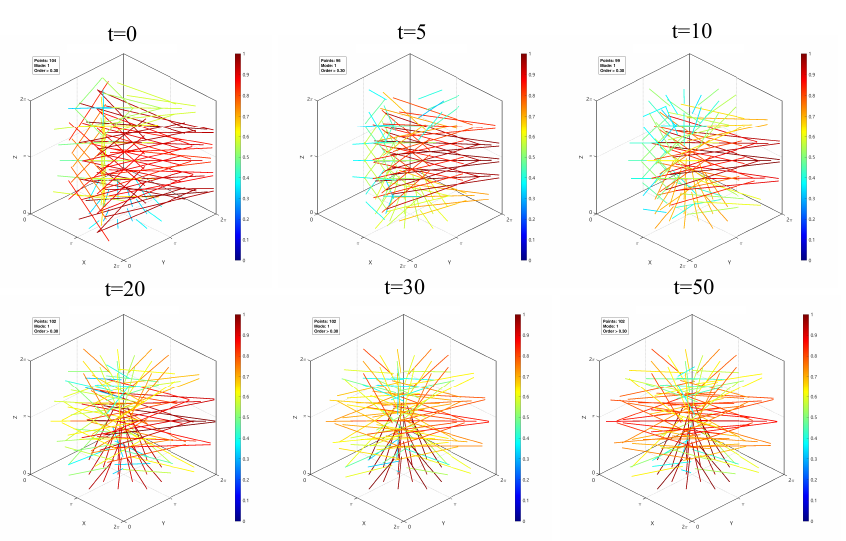}
		 \caption{Evolution trends of the director field  of $Q$ at t = 0,5,10,20,30,50.}
	 \label{v3_deli_qui11}
 \end{figure}
Figure \ref{v3_deli_qui11} illustrates the temporal evolution of the three-dimensional director field, with regions of low scalar order parameter ($S \le 0.30$) omitted for clarity.
Starting from $t=0$, the directors undergo significant reorientation to minimize the system's free energy.
By $t=50$, a stable \textbf{radial hedgehog defect} emerges at the domain center, characterized by a pure splay deformation with directors radiating outward.
The visible ``void'' at the center corresponds to the isotropic defect core ($S \approx 0$), where the nematic order melts to regularize the topological singularity.
 \section{Conclusion}\label{section5}~\\
 In this study, we develop and analyze first- and second-order  low regularity integrators  for  the Q-tensor gradient flow system, which is a fundamental model that describes  liquid crystal  dynamics.
By analyzing the derivative of  tensor functions and their subsequent contraction, we derived two second-order, tensor-based low-regularity integrator (LRI) schemes.
The proposed LRI schemes are both proved to preserve the critical physical constraints of the MBP, energy stability and have optimal convergence rates under temporally continuous regularity assumptions rather than $H^2$ or $H^4$ in other exponential integrator schemes.   Comprehensive numerical experiments validate our theoretical findings, confirming that the proposed LRI schemes effectively preserve physical constraints while maintaining higher accuracy and computational efficiency than exponential time differencing schemes. This work provides a rigorous numerical framework for simulating liquid crystal dynamics with enhanced reliability in preserving essential physical properties.
\section*{Acknowledgments}
G. Ji is partially supported by the National Natural Science Foundation of China (Grant No. 12471363).


\end{document}